\documentclass[12pt]{amsart}
\usepackage{graphicx}
\usepackage{amssymb}
\usepackage{amsfonts}
\usepackage{amsmath}
\usepackage{array}
\usepackage{rotating}

\usepackage{tableaux}
\usepackage[matrix,frame]{xypic}
\newdimen\vcadre\vcadre=0.1cm 
\newdimen\hcadre\hcadre=0.1cm 

\def\arx#1[#2]{\ifcase#1 \relax \or%
  \ar @{-}[#2]  \or%
  \ar @2{-}[#2] \or%
  \ar @{--}[#2] \or%
  \ar @2{.}[#2] \or%
  \ar @{~}[#2]  \fi}

\usepackage{tikz}
\usetikzlibrary{arrows.meta}

\newtheorem{example}{Example}[section]

\newtheorem{theorem}[example]{Theorem}

\newtheorem{corollary}[example]{Corollary}
\newtheorem{conjecture}[example]{Conjecture}
\newtheorem{definition}[example]{Definition}
\newtheorem{proposition}[example]{Proposition}

\newtheorem{lemma}[example]{Lemma}

\def\wt{{\rm wt}}
\def\trg{{\rm tr}}
\def\u{{\bf u}}
\def\wasc{{\rm wasc}}

\def\arm{\operatorname{arm}}
\def\leg{\operatorname{leg}}
\def\NCLLT{{\bf LLT}}
\def\HH{{\bf H}}
\def\Inv{{\rm Inv}}
\def\Asc{{\rm Asc}}
\def\des{{\rm des}}

\def\Proof{\noindent \it Proof -- \rm}
\def\qed{\hspace{3.5mm} \hfill \vbox{\hrule height 3pt depth 2 pt width 2mm}
\bigskip}

\def\XX{{\bf X}}

\def\WQSym{{\bf WQSym}}

\def\<{\langle}
\def\>{\rangle}

\def\C{{\mathbb C}}

\def\N{{\bf N}}
\def\NN{{\mathbb N}}

\def\tr{\operatorname{tr}}

\def\M{{\bf M}}

\def\SG{{\mathfrak S}}

\def\Des{\operatorname{Des}}

\def\maj{{\rm maj}}

\def\GPG{{\mathcal{G}}}
\def\GPD{{\mathcal{D}}}

\def\X{{\rm X}}
\def\LLT{{\rm LLT}}
\def\asc{{\rm asc}}

\def\inv{{\rm inv}}

\makeatletter
\newsavebox{\@brx}
\newcommand{\llangle}[1][]{\savebox{\@brx}{\(\m@th{#1\langle}\)}%
  \mathopen{\copy\@brx\kern-0.5\wd\@brx\usebox{\@brx}}}
\newcommand{\rrangle}[1][]{\savebox{\@brx}{\(\m@th{#1\rangle}\)}%
  \mathclose{\copy\@brx\kern-0.5\wd\@brx\usebox{\@brx}}}
\makeatother

\def\Tabvrule{\vrule width-0.4pt}       
\def\Tabhrule{\hrule \hrule height-0.4pt} 
\def\Tabstrut{\vrule height2.2ex 
                     depth0.8ex  
                     width0ex    
\relax}

\def\PasCase#1{\omit%
            $\vcenter{\hbox {\vbox to 0.4pt{}}
               \hbox{\makebox[3ex]{\Tabstrut$#1$}}}%
               \Tabvrule$}
\def\PasCasePoint{\PasCase{\cdot}}
\def\DessinCarre#1{%
    \vcenter{\hbox{}\hrule
             \hbox{\vrule\makebox[3ex]{\Tabstrut$#1$}\vrule}\Tabhrule}%
             \Tabvrule}
\def\GenRuban#1{\vcenter{\halign{&$\DessinCarre{##}$\cr#1}}\egroup}

\def\sTabvrule{\vrule width-0.4pt}
\def\sTabhrule{\hrule \hrule height-0.4pt}
\def\sTabstrut{\vrule height1.6ex depth0.6ex width0ex \relax}
\def\sDessinCarre#1{%
    \vcenter{\hbox{}\hrule
             \hbox{\vrule\makebox[2.3ex]%
                  {\sTabstrut$\scriptstyle#1$}\vrule}\sTabhrule}%
             \sTabvrule}
\def\sGenRuban#1{\vcenter{\halign{&$\sDessinCarre{##}$\cr#1}}\egroup}

\def\ruban{%
  \bgroup
  \let\ =\omit
  \let\\=\cr
  \let\x=\times
  \let\.=\PasCasePoint
  \offinterlineskip
  \GenRuban}

\def\sruban{%
  \bgroup
  \let\ =\omit
  \let\x=\times
  \let\\=\cr
  \offinterlineskip
  \sGenRuban}

\newdimen\Squaresize \Squaresize=14pt
\newdimen\Thickness \Thickness=0.5pt

\def\Square#1{\hbox{\vrule width \Thickness
   \vbox to \Squaresize{\hrule height \Thickness\vss
      \hbox to \Squaresize{\hss#1\hss}
   \vss\hrule height\Thickness}
\unskip\vrule width \Thickness}
\kern-\Thickness}

\def\Vsquare#1{\vbox{\Square{$#1$}}\kern-\Thickness}

\def\young#1{
\vbox{\smallskip\offinterlineskip
\halign{&\Vsquare{##}\cr #1}}}

\def\boxit#1#2{\setbox1=\hbox{\kern#1{#2}\kern#1}%
\dimen1=\ht1 \advance\dimen1 by #1 \dimen2=\dp1 \advance\dimen2 by #1
\setbox1=\hbox{\vrule height\dimen1 depth\dimen2\box1\vrule}%
\setbox1=\vbox{\hrule\box1\hrule}%
\advance\dimen1 by .4pt \ht1=\dimen1
\advance\dimen2 by .4pt \dp1=\dimen2 \box1\relax}

\def\PC{\rm PC\,}

\def\taille{.5}

\def\gtroistrois{
\begin{tikzpicture}
\begin{scope}[every node/.style={circle,scale=.5,fill=white,draw}]
    \node (A) at (0,0) {};
    \node (B) at (1*\taille,0) {};
    \node (C) at (2*\taille,0) {};
\end{scope}

\begin{scope}[>={Stealth[black]},
              every edge/.style={draw=black,thick}]
    \path [-] (B) edge (C);
\end{scope}
\end{tikzpicture}
}

\def\gtroisquatre{
\begin{tikzpicture}
\begin{scope}[every node/.style={circle,scale=.5,fill=white,draw}]
    \node (A) at (0,0) {};
    \node (B) at (1*\taille,0) {};
    \node (C) at (2*\taille,0) {};
\end{scope}

\begin{scope}[>={Stealth[black]},
              every edge/.style={draw=black,thick}]
    \path [-] (A) edge (B);
    \path [-] (B) edge (C);
\end{scope}
\end{tikzpicture}
}

\def\gquatredeux{
\begin{tikzpicture}
\begin{scope}[every node/.style={circle,scale=.5,fill=white,draw}]
    \node (A) at (0,0) {};
    \node (B) at (1*\taille,0) {};
    \node (C) at (2*\taille,0) {};
    \node (D) at (3*\taille,0) {};
\end{scope}

\begin{scope}[>={Stealth[black]},
              every edge/.style={draw=black,thick}]
    \path [-] (C) edge (D);
\end{scope}
\end{tikzpicture}
}

\def\gquatrecinq{
\begin{tikzpicture}
\begin{scope}[every node/.style={circle,scale=.5,fill=white,draw}]
    \node (A) at (0,0) {};
    \node (B) at (1*\taille,0) {};
    \node (C) at (2*\taille,0) {};
    \node (D) at (3*\taille,0) {};
\end{scope}

\begin{scope}[>={Stealth[black]},
              every edge/.style={draw=black,thick}]
    \path [-] (B) edge (C);
    \path [-] (C) edge (D);
\end{scope}
\end{tikzpicture}
}

\def\gquatresix{
\begin{tikzpicture}
\begin{scope}[every node/.style={circle,scale=.5,fill=white,draw}]
    \node (A) at (0,0) {};
    \node (B) at (1*\taille,0) {};
    \node (C) at (2*\taille,0) {};
    \node (D) at (3*\taille,0) {};
\end{scope}

\begin{scope}[>={Stealth[black]},
              every edge/.style={draw=black,thick}]
    \path [-] (A) edge (B);
    \path [-] (C) edge (D);
\end{scope}
\end{tikzpicture}
}

\def\gquatreneuf{
\begin{tikzpicture}
\begin{scope}[every node/.style={circle,scale=.5,fill=white,draw}]
    \node (A) at (0,0) {};
    \node (B) at (1*\taille,0) {};
    \node (C) at (2*\taille,0) {};
    \node (D) at (3*\taille,0) {};
\end{scope}

\begin{scope}[>={Stealth[black]},
              every edge/.style={draw=black,thick}]
    \path [-] (A) edge (B);
    \path [-] (B) edge (C);
    \path [-] (C) edge (D);
\end{scope}
\end{tikzpicture}
}

\title[Noncommutative chromatic quasi-symmetric functions]{Noncommutative %
chromatic quasi-symmetric functions, Macdonald polynomials, and the %
Yang-Baxter equation}

\author[J.-C.~Novelli and J.-Y.~Thibon]%
{Jean-Christophe Novelli and Jean-Yves Thibon}

\address[] {LIGM, Universit\'e
Gustave-Eiffel, CNRS, ENPC, ESIEE-Paris \\
5 Boulevard Descartes \\Champs-sur-Marne \\77454 Marne-la-Vall\'ee cedex 2 \\
FRANCE}
\email[Jean-Christophe Novelli]{jean-christophe.novelli@univ-eiffel.fr}
\email[Jean-Yves Thibon]{jean-yves.thibon@univ-eiffel.fr} 



\date{}

\begin{document}

\begin{abstract}
As shown in our paper [JCTA  177 (2021), Paper No. 105305], the chromatic
quasi-symmetric function of Shareshian-Wachs can be lifted to $\WQSym$, the
algebra of quasi-symmetric functions in noncommuting variables.  We
investigate here its behaviour with respect to classical transformations of
alphabets and propose a noncommutative analogue of Macdonald polynomials
compatible with a noncommutative version of the Haglund-Wilson formula.  We
also introduce a multi-$t$ version of these noncommutative analogues. For
rectangular partitions, their commutative images at $q=0$ appear to coincide
with the multi-$t$ Hall-Littlewood functions introduced in  [Lett. Math. Phys.
35 (1995), 359--374].  This leads us to conjecture that for rectangular partitions,
multi-$t$ Macdonald polynomials are obtained as equivariant traces of certain
Yang-Baxter elements of Hecke algebras. We also conjecture that all (ordinary)
Macdonald polynomials can be obtained in this way.  We conclude with some
remarks relating various aspects of quasi-symmetric chromatic functions to
calculations in Hecke algebras.
In particular, we show that all modular relations are given by the product
formula of the Kazhdan-Lusztig basis.
\end{abstract}

\maketitle

\section{Introduction}

This paper is a continuation of \cite{NCLLT}. In this reference, we obtained
a noncommutative analogue of the Carlsson-Mellit identity \cite{CM}
relating unicellular LLT polynomials and the
quasi-symmetric chromatic polynomials \cite{SW} of Dyck graphs
\begin{equation}\label{eq:X2LLT}
\X_G(t,X) = (t-1)^{-n}\,\LLT_G(t,(t-1)X).
\end{equation}
We first showed how to deduce this relation from a morphism from the
Guay-Paquet Hopf algebra to $QSym$, and then obtained a noncommutative version
by extending this morphism to $\WQSym$.

This method can be extended to other transformations of alphabets. In Section
\ref{sec:ex}, we describe the image of the noncommutative chromatic
quasi-symmetric function $\XX_G$ by the extension to $\WQSym$ of the
$\omega$-involution of $QSym$, and the image of the noncommutative unicellular
LLT-polynomials $\NCLLT_G(A;t)$ by the transformation $A\mapsto A(q-1)$.
In Section \ref{sec:mac}, we propose a noncommutative lift of the Macdonald
polynomials $H_\mu$ to $\WQSym$. The rationale for this definition is the
Haglund-Wilson formula expressing the $J$-functions in terms of chromatic
quasi-symmetric functions: if we replace in this formula the ordinary $\X_H$
by their noncommutative version, we obtain a multiplicity-free sum of terms
$q^it^m\M_u$ where $i,m$ are  HHL statistics $\inv$ and $\maj$ of the packed
word $u$ interpreted as the row-reading of a filling of the integer partition
$\mu$ associated with $u$.
This definition is then extended to a multi-$t$-analogue. Next, we relate
Macdonald polynomials to the Yang-Baxter bases of Hecke algebras, and propose
a conjectural expression of the Macdonald polynomials as equivariant traces
of intertwiners. We also relate the $e$-positivity conjecture to properties of
the Yang-Baxter basis.
Finally, we conclude by showing that all modular relations derive from the
product formula of the Kazhdan-Lusztig basis of the Hecke algebra.

{\bf Acknowledgements. } This research has been partially supported by the
project CARPLO of the Agence Nationale de la Recherche (ANR-20-CE40-0007).

\section{General setup}

Recall that the Guay-Paquet Hopf algebra $\GPG$ is based on finite simple
undirected graphs with vertices labelled by the integers from 1 to $n=|V(G)|$.
The product is the shifted concatenation: $G\cdot H = G\cup H[n]$ where $H[n]$
is $H$ with labels shifted by the number $n$ of vertices of $G$.

The parameter $t$ arises in the coproduct. If $G$ is a graph on $n$ vertices
and $w\in[r]^n$, regarded as a coloring of $G$, we denote by $G|_w$ the tensor
product $G_1\otimes\cdots\otimes G_r$ of the restrictions of $G$ to vertices
colored $1,2,\ldots,r$. The $r$-fold coproduct is then
\begin{equation}\label{eq:coprod}
\Delta^r G = \sum_{w\in[r]^n} t^{\asc_G(w)} G|_w,
\end{equation}
where $\asc_G(w)$ is the number of edges $(i,j)$ of $G$ such that $w_i<w_j$.
It is also proved in \cite{GP} that the subspace $\GPD$ of $\GPG$ spanned by
Dyck graphs is a Hopf subalgebra.

Given a graph $G$, define the noncommutative chromatic quasi-symmetric
function by
\begin{equation}
\XX_G(t,A) =  \sum_{c\in\PC(G)}t^{\asc_G(c)}\M_{c}(A) \in \WQSym.
\end{equation}
Then, \cite[Prop. 6.1]{NCLLT},
\begin{proposition}
$G\mapsto \XX_G(A)$ is a morphism of Hopf algebras from $\GPG$ to $\WQSym$.
\end{proposition}

As a consequence, for an arbitrary transformation of alphabets (see
\cite[Sec. 3.2]{NCLLT}), the value of $\XX_G(AT)$ where $T$ is a classical
commuting alphabet is determined by the $\XX_H(T)=\X_H(T)$, where $\X_H(X)$ is
the usual quasisymmetric chromatic function.
Indeed, denoting by $\N_u$ the dual basis of $\M_u$,
\begin{equation}
\X_G(T) = \<X_G(X),S_n(TA)\> = \<X_G(XT),S_n(A)\> = \<\XX_G(AT),\N_{1^n}\>
\end{equation}
which is equal to the coefficient of $\M_{1^n}$ in $\XX_G(AT)$. Then, arguing
as in the proof of \cite[Theorem 6.2]{NCLLT}, we can state:

\begin{theorem}\label{th1}
The coefficient of $\M_v$ in $\XX_G(AT)$ is
\begin{equation}
c_v(T) = \<\N_v,\XX_G(AT)\> = t^{\asc_G(v)}\prod_{i=1}^{\max(v)}\X_{G_i(v)}(T).
\end{equation}
\end{theorem}

\begin{lemma}\label{lem1}
Let $u,v$ be packed words and assume that $u$ is finer than $v$. Then,
\begin{equation}
\asc_G(u)= \asc_G(v)+\sum_{i=1}^{\max(v)}\asc_{G_i(v)}(\left(u|_{v^{-1}(i)}\right).
\end{equation}
\end{lemma}
\Proof
By definition of the refinement order, if $v_i>v_j$ then $u_i>u_j$ so that
$\Asc_G(v)\subseteq \Asc_G(u)$,e and $\Asc_G(u)\backslash \Asc_G(v)=\{(i,j)\in
E(G)|v_i=v_j\}$, again because $v_i>v_j$ implies $u_i>u_j$.
\qed

\begin{lemma}
Let $v$ be a packed word. Then, the proper packed colorings $u$ of $G$ 
finer than $v$ are in bijection with sequences $(w_1,\ldots,w_r)$
($r=\max(v)$) of words over $[r]$ such that each $w_i$ is a proper coloring of
$G_i(v)$. The bijection is given by $w_i=u|_{v^{-1}(i)}$.
\end{lemma}

\Proof
Given such a sequence, there is only one $u$ satisfying the required
conditions: the $u_k$ such that $v_k=1$ are given by $w_1$, those such that
$v_k=2$ are obtained by shifting $w_2$ by $\max\left(u|_{v^{-1}(1)}\right)$,
and so on.
\qed

\Proof (of the theorem).

\begin{align}
\XX_G(AT)&=\sum_{u\in\PC(G)}t^{\asc_G(u)}\M_u(AT)\\
&=\sum_{u\in\PC(G)}t^{\asc_G(u)}\sum_{v\le u}\M_v(A) (M_{u|_{v^{-1}(1)}}\cdots M_{u|_{v^{-1}(r)}})(T)\\
&=\sum_v \M_v(A)\sum_{u\ge v} t^{\asc_G(u)}(M_{u|_{v^{-1}(1)}}\cdots M_{u|_{v^{-1}(r)}})(T)\\
&=\sum_v \M_v(A) t^{\asc_G(v)}\sum_{u\ge v} (t^{\asc_{G_1(v)}(w_1)}M_{u|_{v^{-1}(1)}}\cdots t^{\asc_{G_r(v)}(w_r)} M_{u|_{v^{-1}(r)}})(T)\\
& \text{(by Lemma \ref{lem1}, where $w_i=u|_{v^{-1}(i)}$)} \nonumber\\
&=\sum_v\M_v(A) t^{\asc_G(v)}\prod_{i=1}^{\max(v)}\sum_{w\in\PC(G_i(v))}t^{\asc_{G_i(v)}(w)}\M_w(T)\\
&= \sum_v \M_v(A)t^{\asc_G(v)}\prod_{i=1}^{\max(v)}\X_{G_i(v)}(T).
\end{align}
\qed

\section{Examples}\label{sec:ex}

\subsection{The $\omega$-involution}

On symmetric functions, the $\omega$ involution is defined by
$\omega(p_n)=(-1)^{n-1}p_n=(-1)^np_n(-X)$. It can be successively
extended to noncommutative symmetric functions by defining
$\sigma_1(-A)=\sigma_1(A)^{-1}$, so that $S_n(-A)=(-1)^n\Lambda_n(A)$,
then to $QSym$ by duality, to $\WQSym^*$ by using the internal product,
and finally to $\WQSym$ by duality. The result is \cite{NTsuper}
\begin{equation}
\omega \M_u = (-1)^{n-r}\sum_{v\le u}\M_v = (-1)^n \M_u(-A),
\end{equation}
where the sum runs over the refinement order on packed words\footnote{$v\le u$
iff the set composition encoded by $v$ is obtained by merging adjacent blocks
of that encoded by $u$.}, $n=|u|$ and $r=\max(u)$.

{\footnotesize
For example,
\begin{equation}
  \omega \M_{2132}=-(\M_{2132}+\M_{2122}+\M_{1121}+\M_{1111}),
\end{equation}
where the packed words in the r.h.s. encode respectively the set compositions
$2|14|3$, $2|134$, $124|3$, and $1234$.
}

Now, for a Dyck graph $G$, the coefficient of $\M_{1^n}$ in $\omega\XX_G(A)$
is
\begin{equation}
  \<\omega\XX_G,\N_{1^n}\>
 =\<\omega\X_G,S_n\>=\<\omega \X_G,h_n\>
 =\sum_{\lambda\vdash n}\<\omega \X_G,m_\lambda\> =\omega\X_G(1),
\end{equation}
since $\X_G$ is symmetric, that is, the sum of the coefficients of $\omega
\X_G$ on the $h$-basis, which, by \cite[Theorem 5.3]{SW}, is the generating
polynomial of acyclic orientations of $G$ counted by the number of increasing
arcs.

Therefore, we can state: 
\begin{proposition}
The coefficient of $\M_u$ in $\omega\XX_G$ is
\begin{equation}
\<\omega\XX_G,\N_u\>=\sum_{\gamma}t^{\asc_G(\gamma,u)},
\end{equation}
where the sum runs over acyclic orientations $\gamma$ of $G$ without
$u$-descent, that is, without the arcs $(i,j)$ if $u_j<u_i$, and
$\asc_G(\gamma,u)$ is the number of arcs $(i,j)$ of $\gamma$ such that $i<j$.
\end{proposition}

\Proof
By Theorem \ref{th1}, 
\begin{equation}
\<\XX_G(-A),\N_u\> = t^{\asc_G(u)}\prod_{i=1}^{\max(u)}\X_{G_i(u)}(-1)
\end{equation}
so that
\begin{equation}
\<\omega\XX_G(A),\N_u\> = t^{\asc_G(u)}\prod_{i=1}^{\max(u)}\omega\X_{G_i(u)}(1),
\end{equation}
and $\omega\X_H(1)$ is the generating polynomial of acyclic orientations of
$H$ counted by number of increasing arcs.
Thus, $\<\omega\XX_G,\N_u\>$ is the product of $t^{\asc_G(u)}$ by the
generating polynomial of acyclic orientations of the union of the $G_i(u)$ by
number of increasing arcs. This amounts to counting acyclic orientations of
$G$ such that all arcs between different $G_i(u)$ are increasing, that is,
without $u$-descents.
\qed

This is the natural $t$-analogue of a result of Stanley \cite{Stan}.

{\footnotesize
For example, for the graph 
\begin{equation}
	G=
\begin{tikzpicture}
\begin{scope}[every node/.style={circle,scale=.5,fill=white,draw}]
    \node (B) at (1*\taille,0) {};
    \node (C) at (2*\taille,0) {};
    \node (D) at (3*\taille,0) {};
    \node (E) at (4*\taille,0) {};
\end{scope}

\begin{scope}[>={Stealth[black]},
              every edge/.style={draw=black,thick}]
    \path [-] (B) edge (C);
    \path [-] (C) edge (D);
    \path [-] (D) edge (E);
    \path [-] (C) edge[bend left=60] (E);
\end{scope}
\end{tikzpicture}
\end{equation}
labelled $1,2,3,4$ from left to right, the coefficient of $\M_u=\M_{1122}$ in
$\omega\XX_G$ is the product of $\X_{G_1(u)}(1)$ and $\X_{G_2(u)}(1)$, both
equal to $1+t$, by $t^{\asc_G(1122}=t^2$, thus equal to $t^2(1+t)^2$,
corresponding to the four acyclic orientations of $G$ with arcs $(2,3)$ and
$(2,4)$ mandatory, and the others free. Similarly, the coefficient of
$\M_{1212}$ is $t^2+t^3$, with arcs $(1,2),(3,2),(3,4)$ imposed, and the
orientation of $\{2,4\}$ free.
For $u=2211$, $u$ has no $G$-ascent and the result is $(1+t)^2$.
}

\subsection{The $(q-1)$-transform of unicellular LLT polynomials}

Having dealt with the cases $T=-1$ and $T=1/(t-1)$, we can now consider
$T=\frac{q-1}{t-1}$.
\begin{theorem}
For a Dyck graph $G$, let $a_j$ be the number of edges $(i<j)$. Then,
\begin{equation}
(t-1)^n	\X_G\left(\frac{q-1}{t-1}\right) = \prod_{j=1}^n(q-t^{a_j}).
\end{equation}
\end{theorem}

\Proof
Let $m_i(c)$ be the number of occurrences of the letter $i$ in the word $c$,
and denote by $e(G)$ the number of edges of a graph $G$.
Since $G\mapsto \X_G$ is a morphism of Hopf algebras,
\begin{equation}
\X_G\left(\frac{q-1}{t-1}\right)
= \sum_{c\in [2]^n}t^{\asc_G(c)}
    \X_{G_1}\left(\frac{q}{t-1}\right)
    \X_{G_2}\left(\frac{-1}{t-1}\right),
\end{equation}
where, by \cite[Prop. 4.9]{NCLLT}
\begin{equation}
\X_{G_1}\left(\frac{q}{t-1}\right)= q^{m_1(c)}(t-1)^{-m_1(c)}\ \text{and}\
\X_{G_2}\left(\frac{-1}{t-1}\right)=(-1)^{m_2(c)}t^{e(G_2)}(t-1)^{-m_2(c)}
\end{equation}
so that
\begin{equation}
(t-1)^n \X_G\left(\frac{q-1}{t-1}\right)=\sum_{c\in [2]^n}t^{\asc_G(c)}q^{m_1(c)}(-1)^{m_2(c)}t^{e(G_2)}.
\end{equation}
Each monomial in this sum corresponds to one of the $2^n$ terms in the
expansion of the product
\begin{equation}
(q-t^{a_1})(q-t^{a_2})\cdots (q-t^{a_n})
\end{equation}
formed by choosing $q$ in the factors corresponding to the vertices of $G_1$
and $-t^{a_j}$ for those corresponding to the vertices of $G_2$. The power of
$q$ and the sign are obvious, and $t^{\asc_G(c)}$ is the contribution of the
erased edges between $G_1$ and $G_2$.
\qed

{\footnotesize
For example, for the graph 
\begin{equation}
	G=
\begin{tikzpicture}
\begin{scope}[every node/.style={circle,scale=.5,fill=white,draw}]
    \node (B) at (1*\taille,0) {};
    \node (C) at (2*\taille,0) {};
    \node (D) at (3*\taille,0) {};
    \node (E) at (4*\taille,0) {};
\end{scope}

\begin{scope}[>={Stealth[black]},
              every edge/.style={draw=black,thick}]
    \path [-] (B) edge (C);
    \path [-] (B) edge[bend left=60] (D);
    \path [-] (C) edge (D);
    \path [-] (C) edge[bend left=60] (E);
    \path [-] (D) edge (E);
    \path [-] (C) edge[bend left=60] (E);
\end{scope}
\end{tikzpicture}
\end{equation}
the coloring 2121 selects in the product
\begin{equation}
(q-1)(q-t)(q-t^2)(q-t^2)
\end{equation}
the factors $-1,q,-t^2,-t^2$, whose product $-qt^4$ is obtained as
$t^2qt^2(-1)^3$, the first $t^2$ corresponding to the erased edges $(2,3)$ and
$(2,4)$ and the second one to the edges $(1,3)$ and $(3,4)$ of $G_2$.
}

Since this is also the coefficient of $\M_{1^n}$ in 
$(t-1)^n\XX_G(AT)=\NCLLT_G(A(q-1))$, we have:

\begin{corollary}
The coefficient of $\M_u$ in $\NCLLT_G(A(q-1))$ is
\begin{equation*}
t^{\asc_G(u)}\prod_i\LLT_{G_i(u)}\left(q-1\right).
\end{equation*}
\end{corollary}

Setting $q=0$, we obtain:

\begin{corollary}\label{cor:omegaLLT}
The coefficient of $\M_u$ in $\omega\NCLLT_G(A;t)$ is $t^{\wasc_G(u)}$, where
$\wasc_G(u)$ is the number of weak ascents of $u$, i.e., the number of edges
$i<j$ such that $u_i\le u_j$.
\end{corollary}

Indeed, $\<\NCLLT_G(-A),\N_u\>=t^{\asc_G(u)}\prod_i\LLT_{G_i(u)}(-1)$, and 
\begin{equation}
\LLT_H(-1)=(t-1)^{|H|}\X_H\left(\frac{-1}{t-1}\right)=(-1)^{|H|}t^{e(H)}.
\end{equation}
Thus, each $G_i(u)$ contributes by its size to the power of $t$, which amounts
to count not only the $G$-ascents of $u$, but also the edges such that
$u_i=u_j$.

\begin{corollary}
On the $\Phi$-basis of $\WQSym$,
\begin{equation}
\omega\NCLLT_G = \sum_{\sigma\in\SG_n}t^{\wasc_G(\sigma)}\Phi_{\min(\sigma)},
\end{equation}
where the packed word $\min(\sigma)$ is the minimal destandardization
of $\sigma$, {\it e.g.}, $\min(12453)=11221$ (see \cite[Sec. 8]{NCLLT}).
\end{corollary}

\section{Noncommutative Macdonald polynomials}
\label{sec:mac} 

There are many versions of the Macdonald polynomials. We shall first consider
the basis $J$.
In \cite{HW}, Haglund and Wilson define two graphs $G_\mu$ and $G_\mu^+$ 
and prove that
\begin{align}
J_{\mu^{\prime}}(X;q,t) &= t^{-n(\mu^{\prime}) + \binom{\mu_1}{2}} (1-t)^{\mu_1} \sum_{G_{\mu} \subseteq H \subseteq G_{\mu}^{+}} \X_{H}(X;t) \\
&\times \prod_{\{u,v\} \in H \setminus G_{\mu}} - \left(1 - q^{\leg_{\mu}(u) + 1} t^{\arm_{\mu}(u)}\right) \nonumber \\
&\times \prod_{\{u,v\} \in G_{\mu}^{+} \setminus H} \left(1 - q^{\leg_{\mu}(u)+1}t^{\arm_{\mu}(u)+1}\right) . \nonumber
\end{align}
Inserting the Carlsson-Mellit relation
\begin{equation}
\X_G(X;t) = (t-1)^{-n}\LLT_G((t-1)X;t),
\end{equation}
we obtain an expression of $H_\mu$ in terms of $LLT$ polynomials.

In order to obtain a noncommutative generalization of that formula, we shall
start with the HHL formula for $H_\mu(X;q,t)$ in terms of the fundamental
quasisymmetric functions:
\begin{equation}
H_\mu(X;q,t)
= \sum_{\sigma\in\SG_n}
  t^{n(\mu)-\maj(\sigma,\mu)}q^{\inv(\sigma,\mu)}F_{C(\sigma^{-1})},
\end{equation}
where $\maj$ and $\inv$ are defined in \cite{HHL}, and define a noncommutative
analogue by replacing
$F_{C(\sigma^{-1})}$ by $\Phi_{\min(\sigma)}$ in $\WQSym$.

Define
\begin{equation}
\HH_\mu(A;q,t)
= \sum_{\sigma\in\SG_n}
  q^{\Inv(\sigma,\mu)}
  \prod_{u\in\Des(\sigma,\mu)}t^{1+\leg(u)}q^{-\arm(u)}
  \cdot\Phi_{\min(\sigma)} \in\WQSym,
\end{equation}
where we used the legs and the arms of the cells of $\mu$ to compute both
statistics $inv$ and $maj$ (see~\cite{HHL} for more details.)
By definition, the commutative image of $\HH_\mu$ is $H_\mu$.

\medskip
{\footnotesize
For example, one can check that
$\HH_{211}(A;q,t)$
is
\begin{equation}
\label{H211phi}
\begin{split}
\HH_{211}(A;q,t) =& \ \
 t^3  \Phi_{1111} + q t^3  \Phi_{1121} + t \Phi_{1211} + t \Phi_{1212} 
  + q t^3  \Phi_{1221} + q t \Phi_{1321} \\
&+ t^2 \Phi_{2111} + t^2 \Phi_{2112} 
+ q t^2  \Phi_{2121} + t \Phi_{2211} + t^2  \Phi_{2122} + t \Phi_{2212} \\
&  + q t^3  \Phi_{2221} + q t^2  \Phi_{2132} + t \Phi_{2312} + qt\Phi_{2321} 
+ q t^2  \Phi_{3121} + \Phi_{3211} \\
& + \Phi_{3212} + q t^2 \Phi_{3221} 
+ \Phi_{3213} + q t^2  \Phi_{3231} + q t \Phi_{3321} + q \Phi_{4321}.
\end{split}
\end{equation}
}

We shall see that this definition is compatible with the formula of
Haglund-Wilson in the following sense:

\begin{theorem}\label{th:ncH}
Let $\widetilde{\NCLLT}_H(A;q)= q^{e(H)}\NCLLT_H(A;q^{-1})$. Then,
\begin{equation}
\begin{split}
\HH_\mu(A;q,t)=&\frac{q^{-n(\mu')+\binom{\mu_1}{2}}}{(q-1)^{n-\mu_1}}
\sum_{G_\mu\subseteq H\subseteq G_\mu^+}\prod_{u\in H\backslash G_\mu}(1-t^{1+\leg(u)}q^{\arm(u)})\\
	&	\times\prod_{u\in G_\mu^+\backslash H}(t^{1+\leg(u)}q^{\arm(u)+1}-1)\omega\widetilde{\NCLLT}_H(A;q)
\end{split}
\end{equation}
\end{theorem}

\begin{lemma}
For a Dyck graph $G$ of size $n$,
\begin{equation}
\omega\widetilde{\NCLLT}_G(A;q)
 = \sum_{\sigma\in\SG_n}q^{\des_G(\sigma)}\Phi_{\min(\sigma)}.
\end{equation}
\end{lemma}

\Proof
This follows from Corollary \ref{cor:omegaLLT} and \cite[Theorem 8.6]{NCLLT}.
\qed

\Proof (of the Theorem)
Let $R_1$ denote the bottom row of the diagram of $\mu$.
Observing that
\begin{equation}
n(\mu)
= \sum_{u\not\in R_1}(1+\leg(u))\quad\text{and}\quad n(\mu')-\binom{\mu_1}{2}
= \sum_{u\not\in R_1}\arm(u),
\end{equation}
this definition can be rewritten as
\begin{equation}
\HH_\mu(A;q,t)= q^{-n(\mu')+\binom{\mu_1}{2}}\sum_{\sigma\in\SG_n}q^{\Inv(\sigma,\mu)}
\prod_{u\in\Asc(\sigma,\mu)}t^{1+\leg(u)}q^{\arm(u)}\cdot\Phi_{\min(\sigma)}
\end{equation}
where $\Asc(\sigma,\mu)$ is the set of ascents $\{u|u>{\rm South}(u)\}$, where
South(u) indicates the cell below $u$.
Set now 
\begin{equation}
\delta(u)
=
\begin{cases}
1&\text{if $u\in\Asc(\sigma,\mu)$,} \\
0&\text{otherwise,}
\end{cases}
\end{equation}
so that we can now write
\begin{align}
	\HH_\mu(A;q,t)&= q^{-n(\mu')+\binom{\mu_1}{2}}\sum_{\sigma\in\SG_n}q^{\Inv(\sigma,\mu)}
	\prod_{u\not\in R_1}\frac{t^{1+\leg(u)q^{\arm(u)}}(q-q^{\delta(u)}) + q^{\delta(u)}-1}{q-1} \Phi_{\min(\sigma)}\\
	&=\frac{q^{-n(\mu')+\binom{\mu_1}{2}}}{(q-1)^{n-\mu_1}}
	\sum_{\sigma\in\SG_n}q^{\Inv(\sigma,\mu)}
\\& \qquad
	\prod_{u\not\in R_1}\left[(1-t^{1+\leg(u)}q^{\arm(u)})q^{\delta(u)} +
(t^{1+\leg(u)}q^{1+\arm(u)}-1)\right]\Phi_{\min(\sigma)}.
\end{align}
Since $\delta(u)$ can be 0 or 1 for all cells of $\mu\backslash R_1$, we can
expand the product and obtain
\begin{equation}
\begin{split}
\HH_\mu(A;q,t)=&
 \frac{q^{-n(\mu')+\binom{\mu_1}{2}}}{(q-1)^{n-\mu_1}}
 \sum_{\sigma\in\SG_n} q^{\Inv(\sigma,\mu)} \Phi_{\min(\sigma)}
\sum_{U\subseteq \mu\backslash R_1}q^{\delta(u)} \\
&\times	\prod_{u\in U}(1-t^{1+\leg(u)}q^{\arm(u)})\prod_{u\in(\mu\backslash
R_1)\backslash U}(t^{1+\leg(u)}q^{1+\arm(u)}-1).
	\end{split}
\end{equation}
We can now see that if $G_\mu$ and $G_\mu^+$ denote the graphs defined in
\cite{HW}, then, the coordinates of cells of $U$ correspond to the edges of
$G_\mu^+\backslash H$ and those of $(\mu\backslash R_1)\backslash U$ to those
of $H\backslash G_\mu$. Thus,
\begin{equation}
\begin{split}
\HH_\mu(A;q,t)=&\frac{q^{-n(\mu')+\binom{\mu_1}{2}}}{(q-1)^{n-\mu_1}}
\sum_{G_\mu\subseteq H\subseteq G_\mu^+}\prod_{u\in H\backslash G_\mu}(1-t^{1+\leg(u)}q^{\arm(u)})\\
&
\times	\prod_{u\in G_\mu^+\backslash H}(t^{1+\leg(u)}q^{\arm(u)}-1)\sum_{\sigma\in\SG_n}q^{\Inv(\sigma,\mu)+\delta(u)}\Phi_{\min(\sigma)}
\end{split}
\end{equation}
and since $\Inv(\sigma,\mu)+\delta(u)=\des_H(\sigma)$, taking into account
Corollary \ref{cor:omegaLLT}, this can be finally rewritten as
\begin{equation}
\begin{split}
\HH_\mu(A;q,t)=&\frac{q^{-n(\mu')+\binom{\mu_1}{2}}}{(q-1)^{n-\mu_1}}
\sum_{G_\mu\subseteq H\subseteq G_\mu^+}\prod_{u\in H\backslash G_\mu()}(1-t^{1+\leg(u)}q^{\arm(u)})\\
	&\times	\prod_{u\in G_\mu^+\backslash H}(t^{1+\leg(u)}q^{\arm(u)+1}-1)\omega\widetilde{\NCLLT}_H(A;q).
\end{split}
\end{equation}
\qed

For example,
\begin{equation}
\HH_{21}(A;q,t)=
t \Phi_{1 1 1} + \Phi_{1 2 1} + qt\Phi_{2 1 2} + qt \Phi_{2 1 1} + \Phi_{2 2 1} + q \Phi_{3 2 1}
\end{equation}
is equal to
\begin{equation}
(qt-1)\omega\widetilde{\NCLLT}_{\left(\gtroistrois\right)}
+(1-t)\omega\widetilde{\NCLLT}_{\left(\gtroisquatre\right)}
\end{equation}
where
\begin{equation}
\omega\widetilde{\NCLLT}_{\left(\gtroistrois\right)}
 =
 \Phi_{1 1 1} + q \Phi_{1 2 1} + \Phi_{2 1 1} + \Phi_{2 1 2} + q \Phi_{2 2 1} + q \Phi_{3 2 1}
\end{equation}
\begin{equation}
\omega\widetilde{\NCLLT}_{\left(\gtroisquatre\right)}
 = 
 \Phi_{1 1 1} + q \Phi_{1 2 1} + q \Phi_{2 1 1} + q \Phi_{2 1 2} + q \Phi_{2 2 1} + q^2  \Phi_{3 2 1}.
\end{equation}

\section{Multi-$t$-analogues and $R$-matrices}

\subsection{Multi $t$ analogues of the Macdonald polynomials}

Note that the expression of $\HH_\mu$ given in Theorem \ref{th:ncH} allows to
define multi-$t$ Macdonald polynomials by replacing the $t^{1+\leg(u)}$
by new variables $t_{1+\leg(u)}$. Repeating verbatim the above argument shows
that we end up with a sum of terms $a_u(q,t_1,\ldots,t_{n-1})\Phi_u$, where
$a_u$ is a single monomial.

\begin{definition}\label{def:ncH}
Let $\mu$ be a partition of length $r$.
We define
\begin{equation}\label{eq:defHH}
\begin{split}
\HH_\mu(A;q,t_1,\ldots,t_{r-1})=&\frac{q^{-n(\mu')+\binom{\mu_1}{2}}}{(q-1)^{n-\mu_1}}
\sum_{G_\mu\subseteq H\subseteq G_\mu^+}\prod_{u\in H\backslash G_\mu()}(1-t_{1+\leg(u)}q^{\arm(u)})\\
&	\times\prod_{u\in G_\mu^+\backslash H}(t_{1+\leg(u)}q^{\arm(u)+1}-1)\omega\widetilde{\NCLLT}_H(A;q)
\end{split}
\end{equation}
\end{definition}

As explained above, the expansion of these multi-$t$ $\HH_{\mu}$ expand
in the $\Phi$ basis as
\begin{proposition}
Let $\mu$ be a partition of length $r$.
Then
\begin{equation}
\HH_\mu(A;q,{\bf t})=\sum_{\sigma\in\SG_n}q^{\Inv(\sigma,\mu)}\prod_{u\in\Des(\sigma,\mu)}t_{1+\leg(u)}q^{-\arm(u)}\cdot\Phi_{\min(\sigma)}.
\end{equation}
\end{proposition}

Let ${\bf t}$ be the sequence of parameters $(t_1,\dots,t_{r-1},\dots)$.

{\footnotesize	
In the case $\mu = 221$, Formula~\eqref{eq:defHH} reads
\begin{equation}
\begin{split}
	(q\!-\!1)^2 \HH_{221}(A;q,{\bf t}) = &\ \
(qt_2-1)(qt_1-1) \omega\widetilde{\NCLLT}_{\gquatredeux} +
(1-t_2)(qt_1-1) \omega\widetilde{\NCLLT}_{\gquatrecinq} \\
& +
(qt_2-1)(1-t_1) \omega\widetilde{\NCLLT}_{\gquatresix} +
(1-t_2)(1-t_1) \omega\widetilde{\NCLLT}_{\gquatreneuf}.
\end{split}
\end{equation}

We have
\begin{equation}
\begin{split}
\omega\widetilde{\NCLLT}_{\gquatredeux} =&\ \ \ \
   \Phi_{1111} + q \Phi_{1121} + \Phi_{1211} + \Phi_{1212} + q \Phi_{1221} +
   q \Phi_{1321} \\
&+ \Phi_{2111} + \Phi_{2112} + q \Phi_{2121} + \Phi_{2211} + \Phi_{2122} +
   \Phi_{2212} \\
&+ q \Phi_{2221} + q \Phi_{2132} + \Phi_{2312} + q \Phi_{2321} + q
   \Phi_{3121} + \Phi_{3211}  \\
&+ \Phi_{3212} + q \Phi_{3221} + \Phi_{3213} + q \Phi_{3231} + q \Phi_{3321} +
   q \Phi_{4321}
\end{split}
\end{equation}

\begin{equation}
\begin{split}
\omega\widetilde{\NCLLT}_{\gquatrecinq} =& \ \ \ \
   \Phi_{1111} + q \Phi_{1121} + q \Phi_{1211} +
   q \Phi_{1212} + q \Phi_{1221} + q^2  \Phi_{1321} \\
&+ \Phi_{2111} + \Phi_{2112} + q \Phi_{2121} +
   q \Phi_{2211} + \Phi_{2122} + q \Phi_{2212} \\
&+ q \Phi_{2221} + q \Phi_{2132} + q \Phi_{2312} +
   q^2  \Phi_{2321} + q \Phi_{3121} + q \Phi_{3211} \\
&+ q \Phi_{3212} + q \Phi_{3221} + q \Phi_{3213} +
   q \Phi_{3231} + q^2  \Phi_{3321} + q^2  \Phi_{4321}
\end{split}
\end{equation}

\begin{equation}
\begin{split}
\omega\widetilde{\NCLLT}_{\gquatresix} =& \ \ \ \
   \Phi_{1111} + q \Phi_{1121} + \Phi_{1211} +
   \Phi_{1212} + q \Phi_{1221} + q \Phi_{1321} \\
&+ q \Phi_{2111} + q \Phi_{2112} + q^2  \Phi_{2121} +
   \Phi_{2211} + q \Phi_{2122} + \Phi_{2212} \\
&+ q \Phi_{2221} + q^2  \Phi_{2132} + \Phi_{2312} +
   q \Phi_{2321} + q^2  \Phi_{3121} + q \Phi_{3211} \\
&+ q \Phi_{3212} + q^2  \Phi_{3221} + q \Phi_{3213} +
   q^2  \Phi_{3231} + q \Phi_{3321} + q^2  \Phi_{4321}
\end{split}
\end{equation}

\begin{equation}
\begin{split}
\omega\widetilde{\NCLLT}_{\gquatreneuf} =& \ \ \ \
   \Phi_{1111} + q \Phi_{1121} + q \Phi_{1211} +
   q \Phi_{1212} + q \Phi_{1221} + q^2  \Phi_{1321} \\
&+ q \Phi_{2111} + q \Phi_{2112} + q^2  \Phi_{2121} +
   q \Phi_{2211} + q \Phi_{2122} + q \Phi_{2212} \\
&+ q \Phi_{2221} + q^2  \Phi_{2132} + q \Phi_{2312} +
   q^2  \Phi_{2321} + q^2  \Phi_{3121} + q^2  \Phi_{3211} \\
&+ q^2  \Phi_{3212} + q^2  \Phi_{3221} + q^2  \Phi_{3213} +
   q^2  \Phi_{3231} + q^2  \Phi_{3321} + q^3  \Phi_{4321}
\end{split}
\end{equation}

\medskip
One easily checks that the whole sum simplifies into
\begin{equation}
\begin{split}
&\ t_1 t_2 \Phi_{1111} + q t_1 t_2 \Phi_{1121} + t_1 \Phi_{1211} +
   t_1 \Phi_{1212} + q t_1 t_2 \Phi_{1221} + q t_1 \Phi_{1321} \\
&+ t_2 \Phi_{2111} + t_2 \Phi_{2112} + q t_2 \Phi_{2121} +
   t_1 \Phi_{2211} + t_2 \Phi_{2122} + t_1 \Phi_{2212} \\
&+ q t_1 t_2 \Phi_{2221} + q t_2 \Phi_{2132} + t_1 \Phi_{2312} +
   q t_1 \Phi_{2321} + q t_2 \Phi_{3121} + \Phi_{3211} \\
&+ \Phi_{3212} + q t_2 \Phi_{3221} + \Phi_{3213} +
   q t_2 \Phi_{3231} + q t_1 \Phi_{3321} + q \Phi_{4321}.
\end{split}
\end{equation}
which is indeed a multi-$t$ version of Equation~\eqref{H211phi}.
}

\begin{definition}\label{def:ncHt}
Let $\mu$ be a partition of length $r$.
We define
\begin{equation}
\begin{split}
  \widetilde{\HH}_\mu(A;q,t_1,\ldots,t_{r-1})=
 &\frac{q^{-n(\mu')+\binom{\mu_1}{2}}}{(q-1)^{n-\mu_1}}
   \sum_{G_\mu\subseteq H\subseteq G_\mu^+}
      \prod_{u\in H\backslash G_\mu()}(t_{1+\leg(u)}-q^{\arm(u)})\\
 & \times\prod_{u\in G_\mu^+\backslash H}(-t_{1+\leg(u)}+q^{\arm(u)+1})
   \omega\widetilde{\NCLLT}_H(A;q).
\end{split}
\end{equation}
\end{definition}
Clearly, these reduce to the usual $\tilde H_\mu$ for $t_i=t^i$.

\bigskip
{\footnotesize
To compute $\tilde H_{211}(q,t_1,t_2)$, we start from the graph
$G=\gquatredeux$, and we compute
\begin{multline}
[1]\times[(t_1-q)(t_2-q)]\omega\widetilde{\NCLLT}_{\gquatredeux}
+ [-(t_1-1)]\times[(t_2-q)]\omega\widetilde{\NCLLT}_{\gquatresix}\\
+ [-(t_2-1)]\times [(t_1-q)]\omega\widetilde{\NCLLT}_{\gquatrecinq}
+ [-(t_1-1)\times -(t_2-1)]\times [1]\omega\widetilde{\NCLLT}_{\gquatreneuf}
\end {multline}
whose commutative image is
\begin{equation}
\tilde H_{211}(q,t_1,t_2)
= s_4+(q+t_1+t_2)s_{31}+(qt_1+t_2)s_{22}
 + (qt_1+qt_2+t_1t_2)s_{211}+qt_1t_2s_{1111}.
\end{equation}
}

Being linear combinations of LLT polynomials, the commutative images of
$\HH_\mu$ and $\widetilde{\HH}_\mu$ are symmetric functions.

\begin{conjecture}
The coefficients of $H_\mu(q,{\bf t})$ (or $\tilde H_\mu(q,{\bf t})$) on the
Schur basis are polynomials with nonnegative integer coefficients.
\end{conjecture}

Note that the conjecture is true when $\mu$ has at most two rows since in this
case, ${\bf t}$ is reduced to $t_1$, so that $H$ and $\widetilde H$ are
the usual Macdonald functions. It is also true when $\mu=(1^n)$
which is given by an explicit formula
\begin{equation}
\widetilde{\HH}_{1^n}(A;q,{\bf t})
= \sum_{\sigma\in\SG_n}
  \left(\prod_{d\in\Des(\sigma)}t_d\right)\Phi_{\min(\sigma)}.
\end{equation}

Tables of the $\tilde H_\mu$ up to  $n=6$ with at least three rows and two
columns are given in Annex \ref{tables}.

Note that the formula for $\widetilde{\HH}_{1^n}(A;q,{\bf t})$ is reminiscent
to the noncommutative symmetric function
\begin{equation}
K_n(A;t_1,\ldots,t_{n-1}) = \sum_{I\vDash n}\prod_{d\in\Des(I)}t_d\cdot R_I(A)
\end{equation}
whose commutative image is $\tilde H_{1^n}(X;q,t_1,\ldots,t_{n-1})$, and which has been
introduced in \cite[Ex. 3.4]{NCSF2}. It has been observed in \cite{LLT95} that
it was the multi-$t$ analogue of the Hall-Littlewood function $Q'_{1^n}$ given
by the crystal graph, and in \cite{NTT}, it has been used to define a
multi-$(1-t)$-transform on noncommutative symmetric functions.

\subsection{The $q=0$ specialization}

If we set $q=0$ in the definition of $H_\mu(q,{\bf t})$, we obtain a multi-$t$
version of the Hall-Littlewood function $Q'_\mu$, which is itself a
$t$-analogue of the product of complete functions $h_\mu$.

Such multi-$t$ versions have been defined in \cite{LLT95} in the case where
$\mu$ is a rectangular partition. The point of that paper was to describe the
Kostka-Foulkes polynomials $K_{\lambda\mu}(t)$ in terms of crystal graphs.
This was done by averaging a certain statistic on Weyl group orbits, reduced
to a single element for rectangular partitions  $\mu=(k^r)$. In this case, it
seemed natural to replace the terms $\prod t^{id_i(t)}$ by monomials
$\prod t_i^{d_i(t)}$, providing a refinement of the so-called generalized
exponents.

Shortly after the preprint of \cite{LLT95} had been released, other
descriptions of the Kostka-Foulkes polynomials in terms of crystal graphs were
obtained independently by other authors, in particular by Nakayashiki and
Yamada \cite{NY}.

These descriptions involve the energy function on the crystal, which is
defined in terms of the combinatorial $R$-matrix. The combinatorial $R$-matrix
is the limit $q\rightarrow 0$ of the full $R$-matrix, which is defined as an
intertwiner between tensor products of evaluation modules of the quantum
affine algebra $U_q(\widehat{\mathfrak{sl}}_n)$.

\subsection{Connection with the irreducible representations of
$U_q(\widehat{\mathfrak{sl}}_n)$}

The irreducible representations of $U_q({\mathfrak{sl}}_n)$, which  coincide,
as vector spaces, with the irreducible representations $V_\lambda$ of
${\mathfrak{sl}}_n$, can be extended to representations $V_\lambda(z)$ of
$U_q(\widehat{\mathfrak{sl}}_n)$ by choosing a scalar $z\in\C$.

Let $\mu$ be a partition of $N$.
For generic values of the $z_i$, the tensor product
\begin{equation}
W = V_{\mu_1}(z_1)\otimes V_{\mu_2}(z_2)\otimes\cdots\otimes V_{\mu_r}(z_{r})
\end{equation}
is irreducible, and isomorphic to its cyclic shift
\begin{equation}
W'= V_{\mu_2}(z_2)\otimes\cdots\otimes V_{\mu_{r}}(z_{r})\otimes V_{\mu_1}(z_1)
\end{equation}
By Schur's lemma, there is a unique (up to a scalar factor) intertwining
operator between these representations. Let $\Phi$ be such an isomorphism, and
let $\gamma$ be the permutation $(r,1,2,\ldots,r-1)$. Define $g_\gamma:\
W'\rightarrow W$ by
\begin{equation}
g_\gamma(v_1\otimes v_2\otimes \cdots\otimes v_r)
 = v_r\otimes v_1\otimes\cdots\otimes v_{r-1}
\end{equation}
so that $R_\mu:=\Phi\circ g_\gamma$ is now an automorphism of $W$, commuting
with the action of $U_q({\mathfrak{sl}_n})$.
It is therefore a linear combination of the projectors $P_\lambda$ on the
isotypic components of $W$:
\begin{equation}
R_\mu = \sum_{\lambda\vdash N}c_{\lambda\mu}(q;z_1,\ldots,z_{r-1})P_\lambda.
\end{equation}

\subsection{Computing the $c_{\lambda\mu}$}

Alternatively, the coefficients $c_{\lambda\mu}(q;z_1,\ldots,z_{r-1})$ can be
computed by means of the Hecke algebra.
Let $H_N(q)$ and $\hat H_N(q)$ denote respectively the Hecke algebra and the
affine Hecke algebra of type $A_{N-1}$. The finite Hecke algebra $H_N(q)$ is
generated by $T_1,\ldots, T_{N-1}$ satisfying the braid relations with
$(T_i-q)(T_i+1)=0$, and $\hat H_N(q)$ has the extra generators $y_1,\ldots,
y_N$ satisfying $y_iT_j=T_jy_i$ if $j\not=i,i+1$ and $T_iy_iT_i=qy_{i+1}$.

The irreducible modules $S_\lambda$ of $H_N(q)$ can be extended to evaluation
modules $S_\lambda(z)$ of $\hat H_N(q)$ by letting $y_1$ act by multiplication
by $z$. Again, for generic values of the $z_i$, the induction product
\begin{equation}
W = S_{\mu_1}(z_1)\odot S_{\mu_2}(z_2)\odot\cdots\odot S_{\mu_r}(z_{r})
\end{equation}
is irreducible, and isomorphic to its cyclic shift
\begin{equation}
W'= S_{\mu_2}(z_2)\odot\cdots\odot S_{\mu_{r}}(z_{r})\odot S_{\mu_1}(z_1).
\end{equation}
Let $\Phi_\mu$ be an isomorphism $W\rightarrow W'$. 
It is unique up to a scalar factor, and can be interpreted as the right
multiplication by an element $\phi_\mu\in\hat H_N(q)$.

\subsection{Characters of $H_N(q)$}

Recall the rule for computing the characters of $H_N(q)$ \cite{Ram}. Define a
connected element as a product of distinct consecutive generators
$T_iT_{i+1}\cdots T_{i+r-1}$.

The character of such an element does not depend on the order of the factors,
and the character of a product of disjoint (i.e., mutually commuting)
connected elements depends only on the partition $\mu$ recording their
respective number of factors plus one, and its value is
\begin{equation}
\langle f,(q-1)^{-\ell(\mu)}h_\mu((q-1)X))\rangle
\end{equation}
where $f$ is the Frobenius characteristic of the representation under
consideration.

One can always reduce the calculation of the trace of a $T_w$ to that of a
linear combination of such elements, by circularly permuting the factors and
applying the Hecke relations \cite{Ram}.
It is this linear combination which is of interest to us.
We call {\it equivariant trace} of $T_w$ the polynomial obtained by setting
\begin{equation}
\tr(T_iT_{i+1}\cdots T_{i+r-1})=\theta_{r+1}
	\quad\text{where }\\ \theta_i(X;q):=\frac{h_i((q-1)X)}{q-1}
\end{equation}
in such a linear combination.

Then, up to a scalar factor,
\begin{equation}
\tr(\phi_\mu)
 = \sum_{\lambda \vdash N}c_{\lambda\mu}(q;z_1,\ldots,z_{r-1})s_\lambda.
\end{equation}

\subsection{Yang-Baxter elements in the Hecke algebra}

The relevant elements are specializations of the Yang-Baxter basis of
$H_N(q)$. Let us recall its definition \cite{LLT97}. Let $\u=(u_1,\ldots,u_N)$
be a vector of ``spectral parameters''.
For an elementary transposition $s_i=(i,i+1)$, set 
\begin{equation}
	Y_i(x,y) = Y_{s_i}(x,y)=\left(\frac{y}{x}-1\right)T_i+q-1
\end{equation}
and for a permutation $\sigma$ 
\begin{equation}\label{eq:Yw}
	Y_{\sigma s_i}(\u)= Y_\sigma(\u)Y_i(u_{\sigma(i)},u_{\sigma(i+1)}).
\end{equation}
Since the $Y_i$ satisfy the Yang-Baxter equation
\begin{equation}
	Y_i(x,y)Y_{i+1}(x,z)Y_i(y,z)=Y_{i+1}(y,z)Y_i(x,z)Y_{i+1}(x,y),
\end{equation}
Equation \eqref{eq:Yw} allows to compute $Y_\sigma$ from any reduced word of
$\sigma$.

\medskip
{\footnotesize
Let us illustrate this on the smallest possible example. Take $W=S_1(z_1)\odot S_1(z_2)$.
The intertwiner is
\begin{equation}
	Y_{21}(z_1,z_2):= \left(\frac{z_2}{z_1}-1\right) T_1+q-1
\end{equation}
and its equivariant trace is
\begin{align}
	\tr	Y_{21}(z_1,z_2)&= \left(\frac{z_2}{z_1}-1\right) \theta_2+(q-1)\theta_{11}\\
	&= \frac{1}{(q-1)z_1}\left[(z_2-z_1)h_2+z_1 h_{11}\right]((q-1)X)\\
	&= \frac{1}{(q-1)z_1}\left[z_2s_2+z_1 s_{11}\right]((q-1)X)\label{eq:mac11}\\
	&= \frac1{z_1}\left[(qz_2-z_1)s_2+(qz_1-z_2)s_{11}\right].
\end{align}

In this last expression, we can recognize the so-called trigonometric
$R$-matrix for $U_q(\widehat{\mathfrak{sl}}_2)$.
Observe that if we set $z_1=1$ and $z_2=t$, the expression between brackets in
\eqref{eq:mac11} becomes the (not very impressive) Macdonald polynomial
$H_{11}$.

Let us now look at the induced module $W = S_1(1)\odot S_1(t_1) \odot
S_1(t_2)$. The intertwiner is the Yang-Baxter element
\begin{equation}
	Y_{231}(1,t_1,t_2)= ((t_1-1)T_1+q-1)((t_2-1)T_2+q-1)
\end{equation}
whose equivariant trace is
\begin{align}
\tr	Y_{231}(1,t_1,t_2)&=(t_2-1)(t_1-1)\theta_3 +(q-1)(t_1+t_2-2)\theta_{21}+(q-1)^2\theta_{111}\\ 
&= \frac{1}{(q-1)}\left[(t_1-1)(t_2-1)h_3+(t_1+t_2-2)h_{21}+h_{111}     \right]((q-1)X)\\
&= \frac{1}{(q-1)}\left[t_1t_2s_3+(t_1+t_2)s_{21}+s_{111}    \right]((q-1)X)\label{eq:mac111}
\end{align}
and again, we recognize in the bracket of \eqref{eq:mac111} the multi-$t$
Macdonald polynomial $H_{111}$.
Hence, applying the inverse $(1-q)$-transform to \eqref{eq:mac111}, we would
get $\tilde H_{111}(X;q,t_1,t_2)$.
}

Actually, following the steps in the examples above, it is not difficult to
prove that the equivariant trace of the Yang-Baxter element
$Y_{23\cdots n1}(1,t_1,\ldots,t_{n-1})$ is always 
\begin{equation}
(q-1)^{-1}H_{1^n}((q-1)X;q,t_1,\ldots,t_{n-1}).
\end{equation}
In fact, different Yang-Baxter elements give rise to the same Macdonald
polynomial, and here we could as well have used $Y_{n12\cdots
n-1}(1,t_{n-1}/t_{n-2},\ldots,t_{n-1}/t_1,t_{n-1})$.
It is this version which can (apparently) be generalized to rectangular
partitions.

\subsection{Yang-Baxter elements and Macdonald polynomials}

Let $\mu$ be a partition of $n$.
Represent it as its Ferrers diagram (French notation).
We associate a permutation $\sigma_\mu$ and a sequence of parameters $v_\mu$
with $\mu$ as follows.

Number the cells of the Ferrers diagram of $\mu$ from 1 to $n$
starting from the bottom corner, proceeding from left to right and from bottom
to top so as to obtain a standard tableau.  The permutation $\sigma_\mu$ is
obtained by reading the cells of the diagram in the following order:
\begin{itemize}
\item Read by rows from top to bottom and from right to left the values that
have no cell above;
\item Read the remaining values by rows from bottom to top and from
right to left.
\end{itemize}

The parameters $v_\mu$ are obtained by reading the cells in the order of
the filling, each parameter being $q^{{\rm row}(c)-1}t^{{\rm col}(c)-1}$.

\medskip
{\footnotesize
For example, with $\mu=(4,3,1,1)$, the filling is
$$ \tableau{9\\ 8 \\ 5&6&7 \\ 1&2&3&4} $$
the permutation is $\sigma_\mu =976432158$, and the parameters corresponding
to the cells are
$$ \tableau{t^3\\ t^2 \\ t&qt&q^2t \\ 1&q&q^2&q^3} $$
so that
$v_\mu=(1,q,q^2,q^3,t,qt,q^2t,t^2,t^3)$.
}

\medskip
\begin{conjecture}
The equivariant trace of the Yang-Baxter element
$Y_{\sigma_\mu}(v_\mu)$ normalized by the coefficient of $s_n$ is equal
to the Macdonald polynomial $\tilde H_\mu((1-q)X;q;t)$.
\end{conjecture}

We have checked this conjecture up to all partitions of size $8$. Actually,
there are many other permutations giving the same Macdonald polynomial.
A subset of those permutations is obtained by shuffling in all possible ways
the rows of the values with no cell above in the definition of $\sigma_{\mu}$
before reading the remaining values in the tableau.

\medskip
{\footnotesize
Here is, for each partition $\mu$ of $6$, the number of permutations yielding
$\tilde H_\mu$ with the parameters $v_\mu$:
$$
(6): 89,\ (51): 40,\ (42): 30,\ (411): 12,\ (33): 21,\ (321): 12,$$
$$(3111): 4,\  (222): 14,\ (2211): 4, \ (21111): 2,\ (111111): 2
$$
and the 12 permutations giving $\tilde H_{321}$ are
$$
 356124,
 356214,
 365124,
 365214,
 536124,
 536214,$$
 $$
 563124,
 563214,
 635124,
 635214,
 653124,
 653214.
$$
}

\medskip
Moreover, one can also ask whether one can recover the multi-$t$ version of
the Macdonald polynomials $\tilde H$ in this way. We were able to find all
multi-t $\tilde H$ as traces of Yang-Baxter elements for partitions up to size
$7$ except for $\tilde H_{2221}$.

\medskip
{\footnotesize
For example, $\tilde H_{211}((1-q)X;q,t_1,t_2)$ is the normalized trace of 
\begin{equation}
\begin{split}
Y_{4213}(1,q,\frac{t_2}{t_1},t_2)
=& \left((t_1-1)T_3+q-1\right)
  \cdot\left(\left(\frac{t_2}{q}-1\right)T_2+q-1\right) \\
&  \cdot\left((t_2-1)T_1+q-1\right)
   \cdot\left((q-1)T_2+q-1\right)
\end{split}
\end{equation}
}

\bigskip
When the partition is a rectangle, we propose the following
conjecture\footnote{This conjecture was formulated in 2001, in collaboration
with A. Lascoux and B. Leclerc, during the semester on Macdonald polynomials
held at the Newton institute. It was part of an unsuccessful attempt to define
multi-$q$ and multi-$t$ Macdonald polynomials (which were empirically
constructed up to $n=5$), and was never published, although mentioned in the
talk \cite{Th}.}

\begin{conjecture}
Let $\mu=(k^r):=(k,k,\ldots,k)$ ($r$ times) be a rectangular partition.
The equivariant trace of the Yang-Baxter element $Y_{\sigma_\mu}(v_\mu)$ 
is proportional to 
\begin{equation}
	\omega	\tilde H_\mu((1-q)X;q,t_1,\ldots,t_{r-1})
\end{equation}
where the vector of spectral parameters is 
\begin{equation}
v_\mu=\left(1,q,\ldots,q^{k-1},  \frac{t_{r-1}}{t_{r-2}},\ldots,q^{k-1}\frac{t_{r-1}}{t_{r-2}},\frac{t_{r-1}}{t_{r-3}},\ldots,q^{k-1}\frac{t_{r-1}}{t_{r-3}},\ldots,t_{r-1},\ldots,q^{k-1}t_{r-1}\right).
\end{equation}
\end{conjecture}
\medskip

{\footnotesize
For example, applying the inverse $(1-q)$ transform to $Y_{652143}\left(1,q,\frac{t_2}{t_1},q\frac{t_2}{t_1},t_2,qt_2\right)$ and dividing by the coefficient of $s_6$,
we obtain
\begin{align*}
\tilde H_{222}=&s_{6} +(q t_1+q t_2+q+t_1+t_2) s_{51}\\
&+(q^2 t_1+q^2 t_2+q t_1 t_2+q^2+q t_1+q t_2+t_1^2+t_1 t_2+t_2^2) s_{42}\\
&+(q^2 t_1 t_2+q^2 t_1+q^2 t_2+q t_1^2+2 q t_1 t_2+q t_2^2+q t_1+q t_2+t_1 t_2) s_{411}\\
	&+(q^3+q t_1^2+q t_1 t_2+q t_2^2+t_1 t_2) s_{33}\\
	&+(q^3 t_1+q^3 t_2+ q^2 t_1^2+2 q^2 t_1 t_2+q^2 t_2^2+q t_1^2 t_2+q t_1 t_2^2+q^2 t_1+q^2 t_2+q t_1^2+2 q t_1 t_2+q t_2^2+t_1^2 t_2+t_1 t_2^2) s_{321}\\
	&+(q^3 t_1 t_2+q^2 t_1^2 t_2+q^2 t_1 t_2^2+q^2 t_1^2+2 q^2 t_1 t_2+q^2 t_2^2+q t_1^2 t_2+q t_1 t_2^2+q t_1 t_2) s_{3111}\\
	&+(q^3 t_1 t_2+q^2 t_1^2+q^2 t_1 t_2+q^2 t_2^2+t_1^2 t_2^2) s_{222}\\
	&+(q^3 t_1^2+q^3 t_1 t_2 +q^3 t_2^2+q^2 t_1^2 t_2+q^2 t_1 t_2^2+q t_1^2 t_2^2+q^2 t_1 t_2+q t_1^2 t_2+q t_1 t_2^2) s_{2211}\\
	&+(q^3 t_1^2 t_2+q^3 t_1 t_2^2+q^2 t_1^2 t_2^2+q^2 t_1^2 t_2+q^2 t_1 t_2^2) s_{21111}+t_1^2 t_2^2 q^3 s_{111111}
\end{align*}
}

It seems also possible to recover multi $t$ Macdonald polynomials for hook
partitions and for two-row partitions, using the circular shift sending the
last factor to the first position.

{\footnotesize
For example, $\tilde H_{211}$ is obtained from $Y_{4213}(1,q,t_2/t_1,t_2)$,
and $\tilde H_{32}$ is obtained from $Y_{54321}(1,q,q^2,t,qt)$.
}

\medskip
However, this does not work for two-column partitions. 

\medskip
{\footnotesize
For example, with $\mu=221$, $Y_{52143}(1,q,t_2/t_1,qt_2/t_1,t_2,qt_2)$
yields a Schur positive polynomial which does not reduce to a Macdonald
function for $t_i=t^i$:
\begin{align*}
s_{5}
	&+(2q+t_1+t_2)s_{41}
+(q^2+qt_1+qt_2+t_1+t_2)s_{32}
+(t_1q^2+q^2t_2+qt_1+qt_2+t_1t_2)s_{221}\\
	&+(q^2+2qt_1+2qt_2+t_1t_2)s_{311}
+q(qt_1+qt_2+2t_1t_2)s_{2111}
+q^2t_1t_2s_{11111}
\end{align*}

But with $Y_{43521}(1,q,t_1,qt_1,t_2)$, we obtain a multi-$t$ analogue of
$\tilde H_{221}$:
\begin{align*}
s_{5} &	+(qt_1+q+t_1+t_2)s_{41}
 + (t_1q^2+qt_1t_2+t_1^2q+qt_1+qt_2+t_1t_2)s_{311} 
	+(q^2+qt_1+qt_2+t_1t_2+t_1^2)s_{32}\\
	&+(t_1q^2+q^2t_2+qt_1t_2+ t_1^2q+t_1^2t_2)s_{221}
	+(qt_1t_2+q^2t_1t_2+t_1^2qt_2+t_1^2q^2)s_{2111}
	+t_1^2q^2t_2s_{11111}
\end{align*}
which however does not coincide with the one given by the Haglund-Wilson
formula: the difference is $(q-1)(t_1^2-t_2)(s_{32}+qs_{221})$.
}

\bigskip
It appears that the inverse $(1-q)$-transform of equivariant traces of
Yang-Baxter elements can produce many Schur-positive symmetric functions,
which share many properties with the family of Macdonald polynomials, though
not always belonging to it. In particular, it seems that the following
conjectural property of our multi-$t$ Macdonald polynomials is also satisfied
by the normalized traces whose coefficients are in $\NN[q,{\bf t}]$.

\begin{conjecture}
Let $\tilde H_\mu(X;q,{\bf t})=\sum_\lambda \tilde K[q,{\bf t}]s_\lambda$ and
write
\begin{equation}
\tilde K[q,{\bf t}]= \sum_{k,\alpha} c_{k,\alpha} q^k {\bf t}^\alpha
\end{equation}
as a sum of monomials. Then,
\begin{equation}
\prod_{k,\alpha}(q^kT^\alpha)^{c_{k,\alpha}}
= \left(\prod_{(i,j)\in\mu} q^{i-1}t_{j-1}\right)^{d_\mu}
\end{equation}
where $(i,j)$ are the cells of the Ferrers diagram of $\mu$ and $d_\mu$ is the
number of standard tableaux of shape $\mu$ in which 2 is above 1, and we set
$t_0=1$.
\end{conjecture}

{\footnotesize
For example, the coefficient of $s_{321}$ in $\tilde H_{222}$ is
$$ (q^3 t_1+q^3 t_2+ q^2 t_1^2+2 q^2 t_1 t_2+q^2 t_2^2+q t_1^2 t_2
   +q t_1 t_2^2+q^2 t_1+q^2 t_2+q t_1^2+2 q t_1 t_2+q t_2^2+t_1^2 t_2+t_1 t_2^2) $$
and the above product yields $(q^3t_1^2t_2^2)^8$.

One can check that the property is also satisfied by the two non-Macdonald
examples associated with $\mu=221$ above.
}

\medskip
Note that even in the case of the usual Macdonald polynomials, this property
seems to hold. We could not find it written in the literature.

\medskip
Note that the coefficient of a hook Schur function $s_{n-k,1^k}$ in $\tilde
H_\mu$ appears to be the elementary symmetric function $e_k(v_\mu-1)$, a known
property of the single-$t$ version.

\section{Generic traces of Bruhat intervals in the Hecke algebra}

\subsection{Traces of Bruhat intervals}

According to \cite{CHSS}, the chromatic quasi-symmetric function of a Dyck
graph $G$ is related  to the equivariant trace of a Bruhat interval in the
Hecke algebra by
\begin{equation}
\omega X_G(q)=  \trg\left(\sum_{v\le w}T_v\right)
\end{equation}
where $w$ is a $312$-avoiding permutation, and the entry $c_i$ of the Lehmer
code $c$ of $w$ is the number of edges $i<j$  in $G$.

Recall that the permutations of the Bruhat interval $[{\rm Id},w]$ have as
reduced words all the reduced subwords of a reduced word of $w$. For a reduced
word $u$, denote by $T(u)$ the corresponding product of generators of the
Hecke algebra, and write $u\sim v$ if $T(u)$ and $T(v)$ have the same trace.

\medskip
{\footnotesize
Here is an example of a direct calculation.
For $w = 3241$, a reduced word is 2123 and the reduced words of the Bruhat
interval, together with their contributions to the trace, are
$$\epsilon \rightarrow \theta_{1111},\ 1,2,3 \rightarrow 3\theta_{211},\ 21, 12, 23 \rightarrow 3\theta_{31},\  13 = 1|3 \rightarrow \theta_{22},\
213,123 \rightarrow 2\theta_{4}$$
$$212 \sim 221 : T(221)=((q-1)T(2)+q)T(1) \rightarrow(q-1)\theta_{31}+q\theta_{211},$$
$$2123 = 1213 = 1231 \sim 1123 : T(1123)= ((q-1)T_1+q)T(23) \rightarrow(q-1)\theta_{4}+q\theta_{31}$$
and finally
$$\trg(T_{3241})= \theta_1^4+(q+3)\theta_{211}+\theta_{22}+(2q+2)\theta_{31}+(q+1)\theta_4.$$
Substituting $t_i=h_i((q-1)X)/(q-1)$ and converting to the $h$-basis, we obtain
$$\left(q^{3} + 2 q^{2} + q\right) h_{3,1} + \left(q^{4} + 2 q^{3} + 2 q^{2} + 2 q + 1\right) h_{4}$$
which is indeed $\omega X_G$, where $G$ is labelled by the partition
$\mu=(2)$, that is
%
$
        G=
\begin{tikzpicture}
\begin{scope}[every node/.style={circle,scale=.5,fill=white,draw}]
    \node (B) at (1*\taille,0) {};
    \node (C) at (2*\taille,0) {};
    \node (D) at (3*\taille,0) {};
    \node (E) at (4*\taille,0) {};
\end{scope}

\begin{scope}[>={Stealth[black]},
              every edge/.style={draw=black,thick}]
    \path [-] (B) edge (C);
    \path [-] (C) edge (D);
    \path [-] (D) edge (E);
    \path [-] (B) edge[bend left=60] (D);
\end{scope}
\end{tikzpicture}
$.
}

\bigskip
Permutations avoiding $312$ are always nonsingular, and the sum over the
Bruhat interval can be factorized as a product of terms $T_i+f(q)$ by an
algorithm due to Lascoux \cite{Las}, which gives a faster way to compute the
equivariant trace directly on the $h$-basis.

\subsection{Factorizations of Lascoux and the Yang-Baxter basis}

Let
\begin{equation}
c_w = \sum_{v\le w}T_w.
\end{equation}
The algorithm of Lascoux \cite{Las} works as follows. A permutation $w$ is nonsingular
iff $w$ or $w^{-1}$ has the property that there exists $k$ such that
$w(k)=n>w(k+1)>\cdots >w(n)$, and $w\backslash n$ is nonsingular.
In the first case, $c_w=c_{w\backslash n}(T_{n-1}+\frac{1}{[n-k}])\cdots
(T_k+\frac{1}{[1]})$.
In the second case,
$c_w = (T_k+\frac{1}{[1]})\cdots (T_{n-1}°\frac1{[n-k]})
       c_{(w^{-1}\backslash n)^{-1}}$.

For 312-avoiding permutations, we are always in the first case.

\medskip
{\footnotesize
Continuing the previous example, we have 
$$c_{3241}=(T_1+1)\left(T_2+\frac1{[2]}\right)(T_1+1)(T_3+1).$$
}

Other factorizations can be deduced from this one by iterated application of
the Yang-Baxter relations.

\subsection{Direct expansion on the $h$-basis}

Let
\begin{equation}
\Upsilon_n
= (T_1+1)\left(T_2+\frac1{[2]}\right)\cdots \left(T_{n-1}+\frac1{[n-1]}\right).
\end{equation}
Then \cite{Las06},
\begin{equation}
\tr \Upsilon_n=[n]_q h_n.
\end{equation}
One can also start with $T_{n-1}$, since the involution $T_i\mapsto T_{n-i}$
is an automorphism of $H_n(q)$.

Indeed, $y_n:=[n-1]!\tr (\Upsilon_n)$ satisfies the recursion
\begin{equation}
y_n = y_{n-1}\theta_1+[n-1]y_{n-2}\theta_2+[n-1][n-2]y_{n-3}\theta_3+\cdots+[n-1]!\theta_n
\end{equation}
which can be seen by grouping the terms of the product
$[n-2]!\Upsilon_{n-1}(1+[n-1]T_{n-1})$ according to the longest connected
right factor $T_iT_{i+1}\cdots T_{n-1}$. 
Thus, if we assume 
$ y_i=[i]!h_i$ for $i<n$, then
\begin{equation}
y_n = [n-1]![t^n]\sigma_t(X)\frac{\sigma_t((q-1)X)-1}{q-1}
=[n-1]![t^n]\frac{\sigma_t(qX)-\sigma_t(X)}{q-1}=[n]!h_n(X).
\end{equation}

In some cases, we can obtain $\omega X_G$ directly in the $h$ basis by
manipulating the factorization so as to transform it into a linear combination
of direct products of factors $\Upsilon_i$.

{\footnotesize
For example,
$$c_{3241}=(T_1+1)\left(T_2+\frac1{[2]}\right)(T_1+1)(T_3+1)\sim (1+q)(T_1+1)\left(T_2+\frac1{[2]}\right)(T_3+1)$$
by moving the second $T_1+1$ to the right and then cycling it. Next, writing the last factor as
$$T_3+1= T_3+\frac1{[3]}+\frac{q+q^2}{[3]},$$
we obtain
$$c_{3241}\sim (1+q)\Upsilon_4+ \frac{q[2][2]}{[3]}\Upsilon_{31}$$
whose trace is
$$[2][4]h_4+q[2][2]h_{31}=\omega X_G.$$
}

\subsection{The formula of Abreu and Nigro}

The main result of \cite{AN} (formulated in terms of chromatic quasi-symmetric
functions) amounts to a combinatorial formula for the generic trace of $c_w$,
where $w$ is a $312$-avoiding permutation. It reads (Th. 1.2. of \cite{AN})
\begin{equation}	
\trg c_w = \sum_{v\le h} q^{\wt_h(v)}\theta_{\lambda(v)}
\end{equation}
where $\lambda(v)$ is the cycle type of $v$, $h$ is the Hessenberg function
associated with $w$, and $v\le h$ means $v(i)\le h(i)$ for all $i$. The
statistic $\wt_h(v)$ is the number of $h$-inversions of the forget-cycles
transform $v'$ of $v$, \emph{i.e.}, pairs $j>i$ such that $v'(i)<h(v'(j))$. The
forget cycles transform is obtained by writing the cycles of $v$ starting with
their smallest element and ordered by increasing values of this element, and
then erasing the parentheses.

The Hessenberg function $h$ is given by $h(i) = c_i+i$, $i=1\ldots n$ where
$c$ is the Lehmer code of $w$, i.e., $c_i$ is the number of $j<i$ to the right
of $i$.

\medskip
{\footnotesize For example the code of $w=3241$ is $c=2110$ so that $h=3344$.
The permutations $v\le_h w$ are with their cycle type and $h$-inversions 
$$
\begin{matrix}
	v&\lambda(v)&v'&{\rm hinv}_(v)&q^{\wt_h(v)}\\
1 2 3 4 & 1 1 1 1 & 1 2 3 4 &  & 1 \\
1 2 4 3 & 2 1 1 & 1 2 3 4 &  & 1 \\
1 3 2 4 & 2 1 1 & 1 2 3 4 &  & 1 \\
1 3 4 2 & 3 1 & 1 2 3 4 &  & 1 \\
2 1 3 4 & 2 1 1 & 1 2 3 4 &  & 1 \\
2 1 4 3 & 2 2 & 1 2 3 4 &  & 1 \\
2 3 1 4 & 3 1 & 1 2 3 4 &  & 1 \\
2 3 4 1 & 4 & 1 2 3 4 &  & 1 \\
3 1 2 4 & 3 1 & 1 3 2 4 & (2 3) & q \\
3 1 4 2 & 4 & 1 3 4 2 & (2 4) & q \\
3 2 1 4 & 2 1 1 & 1 3 2 4 & (2 3) & q \\
3 2 4 1 & 3 1 & 1 3 4 2 & (2 4) & q 
\end{matrix}
$$
which gives back
$$\trg(T_{3241})= \theta_1^4+(q+3) \theta_{211}+\theta_{22}+(2q+2)\theta_{31}+(q+1)\theta_4.$$
}

\subsection{Modular law, linear relations, and the Kazhdan-Lusztig basis}

There are linear relations between certain chromatic quasi-symmetric functions.
One way to formulate them (Guay-Paquet \cite{GP,AN0}) is by viewing  the boundary
of the partition $\mu$ as a Dyck path with east and north steps and  applying when allowed
({\it i.e.}, to an abelian subpath)
the  local tranformations
 $$(1 + q)ene = qeen + nee,$$
 $$(1 + q)nen = qenn + nne.$$
Since we know that $\omega X_G$ is the generic trace of a Bruhat interval in
the Hecke algebra, it is likely that these relations are consequences of the
Hecke relations. 

This fact, which is a straightforward consequence of the properties of
the Kazhdan-Lusztig basis, has been first observed by Abreu and Nigro \cite[Cor. 3.2]{AN1}.
We provide below the precise correspondence with the relations of Lee \cite{Lee}.

When $w$ is $312$-avoiding, and more generally when $w$ is nonsingular, the
interval $c_w$ is, up to a power of $q$, a Kazhdan-Lusztig basis element
$$C'_w(q)=q^{-\ell(w)/2}\sum_{v\le w}P_{v,w}(q)T_v,$$
that is, 
$$c_w=q^{\ell(w)/2}C'_w.$$
The $h$-positivity of the generic trace of all the
$C'_w(q)$ is an old conjecture of Haiman \cite{Hai}.

The basis $C'_w$ satisfies the following relations. If
$s=s_i$ is an elementary transposition,
$$C'_sC'_w=\begin{cases}C'_{sw}+\sum_{sv<v}\mu(v,w)C'_v&\text{if $sw>w$}\\
(q^{1/2}+q^{-1/2})C'_w&\text{if $sw<w$}\end{cases}$$
$$C'_wC'_s=\begin{cases}C'_{ws}+\sum_{vs<v}\mu(v,w)C'_v&\text{if $ws>w$}\\
(q^{1/2}+q^{-1/2})C'_w&\text{if $ws<w$}\end{cases}$$
where $\mu(v,w)$ is the coefficient of $q^{(\ell(w)-\ell(v)-1)/2}$ in $P_{v,w}(q)$ if
$\ell(w)-\ell(v)$ is odd, and $0$ otherwise.

\medskip
{\footnotesize
Here is an example of a modular relation. Let
$w_1= 3 4 2 6 5 1$,  $w_0=3 2 4 6 5 1$, and $w_2= 4 3 2 6 5 1$.
Then, we have the relation
$$(q+1)\trg(c_{w_1})=\trg(c_{w_2})+q\trg(c_{w_0})$$
which can be explained as follows.
The reduced words given by Lascoux's factorization for for $w_1,w_0,w_2$ are
respectively
$$1 2 1 3 2 4 5 4, 1 2 1 3 4 5 4, 1 2 1 3 2 1 4 5 4$$
We can compute
$$c_{w_1}c_{s_1}=c_{w_2}+qc_{w_0}$$
and next, $c_{w_1}c_{s_1}$ is conjugate to $c_{s_1}c_{w_1}=(q+1)c_{w_1}$
which has therefore the same trace.

In this example, $C'_{s_1}C'_{w_1}$ is computed by listing all permutations
$v<w$ such that $vs_1<v$.
This condition is satisfied by $w_0$, and $\ell(w_1)-\ell(w_0)=1$ is odd.
Since $w_1$ is nonsingular, $P_{w_0,w_1}(q)=1$.
The other permutations $v$ 
satisfying these conditions are all strictly less than $w_0$,
since the reduced decomposition of 
$w_1=12132454$ 
forces the removal of the rightmost 2 so as to have
$vs_1<v$, so that for other 
$v$, $\ell(w_1)-\ell(v)\ge 3$ and then
$\mu(v,w_1)=0$ since $P_{v,w_1}(q)=1$.
}

\medskip
\begin{theorem}
All the linear relations described by Lee \cite{Lee}
are of this type.
\end{theorem}

\Proof Recall that Dyck graphs are in bijection with integer partitions included
in the staircase partition as follows: the edges of the graph $G$ are the
coordinates of the empty cells above the diagonal.

{\footnotesize
For example, the partition $(2,2,1)$ represented as
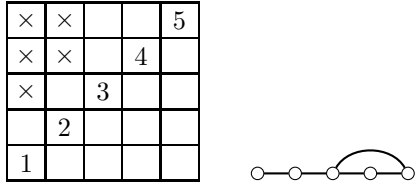
\begin{equation*}
\young{\times &\times & & & 5\cr
       \times & \times & &4& \cr
       \times & & 3 &&\cr
        & 2 &&&\cr
        1&&&&\cr}
\quad\quad
\begin{tikzpicture}
\begin{scope}[every node/.style={circle,scale=.5,fill=white,draw}]
    \node (A) at (0,0) {};
    \node (B) at (1*\taille,0) {};
    \node (C) at (2*\taille,0) {};
    \node (D) at (3*\taille,0) {};
    \node (E) at (4*\taille,0) {};
\end{scope}

\begin{scope}[>={Stealth[black]},
              every edge/.style={draw=black,thick}]
    \path [-] (A) edge (B);
    \path [-] (B) edge (C);
    \path [-] (C) edge (D);
    \path [-] (D) edge (E);
    \path [-] (C) edge[bend left=60] (E);
\end{scope}
\end{tikzpicture}
\end{equation*}
whic corresponds to the graph whose edges are $(1,2),(2,3),(3,4),(3,5),(4,5)$.
}

\medskip
To relate the recursion of Lee to computations in the Hecke algebra, we
will use a slight modification of the presentation given in  \cite{Ale}, by taking the conjugate
partitions and consider column area vectors instead of row area vectors.

We shall encode a partition by its \emph{code sequence}: given $\lambda$
contained in a staircase, its values are the number of empty cells above the diagonal in each
column, read from left to right.
Note that a code sequence is  always the code of a permutation $\sigma$.
In the previous example, the code sequence is $(1,1,2,1,0)$ and $\sigma=23541$.
The condition that an edge $(i,j)$ be admissible ({\it cf.} \cite[Def.  16]{Ale}) for a given code sequence $a=(a_1,\dots,a_k)$
translates now as
\begin{itemize}
\item $j=i+a_i$,
\item $i=1$ or $a_{i-1}\leq a_i-1$,
\item $a_i\geq2$,
\item $a_{j-1} = 1 + a_j$.
\end{itemize}

The permutations $\sigma$ corresponding to these code sequences
are $312$-avoiding, and the inclusion order of the partition diagrams
is the restriction of the Bruhat order to $312$-avoiding permutations.
Indeed, the condition $a_i\geq a_{i-1}+1$
is known (and easily seen) to be the condition satisfied by the codes of permutations avoiding
the pattern $312$.

\begin{lemma}
Consider a code sequence $a$ with an admissible edge $(i,j)$, and let
$\sigma$ be the permutation of code $a$. Let also $\sigma'=s_{j-1}\sigma$.
Then
\begin{itemize}
\item $\sigma_i=j$,
\item $j-1$ is a descent of $\sigma$,
\item $j-1$ is a recoil of $\sigma$,
\item $\sigma'$ also avoids the pattern $312$,
\item $j-1$ is a descent but not a recoil of $\sigma'$.
\end{itemize}
\end{lemma}

\Proof
Since $a_{i-1}<a_i$, we have $\sigma_{i-1}<\sigma_i$. Since $\sigma$ avoids
312, all values $\sigma_\ell$ with $\ell<i$  satisfy
$\sigma_\ell<\sigma_i$.
So in $\sigma$, there are $i-1$ values smaller than $\sigma_i$ on its left and
$a_i$ such values on its left. So $\sigma_i=i+a_i=j$, whence the first item.

Now, since $a_{j-1}=1+a_j$,  $j-1$ is a descent of $\sigma$.
Moreover, all $\ell<i$ satisfy $a_\ell\leq a_{i-1}+i-1-\ell$, so that
$\sigma_\ell\leq a_\ell+\ell \leq a_{i-1}+i-1$.
Since either $i=1$ or $a_{i-1}< a_i$, we have $\sigma_\ell<j-1$, so that
the value $j-1$ is to the right of $j$ in $\sigma$, which means that
$j-1$ a recoil of $\sigma$. 

Define now $\sigma'=s_{j-1}\sigma$ and let $a'$ be its code.
Since we exchanged the values $j-1$ and $j$ in $\sigma$, all values of the
code of $\sigma'$ are the same as those of $\sigma$ (since they have the
same comparison with the values to their right) except for the $i$-th
component  (since $\sigma_i=j$) which is replaced by $a'_i=a_i-1$.
This new sequence also satisfies the condition of being a code sequence, so
$\sigma'$ also avoids the pattern $312$.

Now, since $j-1$ was a recoil of $\sigma$ and since we exchanged $j-1$ and $j$
to get $\sigma'$, then $j-1$ is not a recoil of $\sigma$. Since both values
were not in positions $j-1$ and $j$ in $\sigma$ since $a_i=j$ and
$j-i=a_i\geq2$, there is still a descent in $j-1$ in $\sigma'$. \qed


We can now apply the
product formulas of the Kazhdan-Lusztig basis 
to this situation.

\begin{lemma}
Consider a code sequence $a$ with an admissible edge $(i,j)$, and let
$\sigma$ be the permutation of code $a$ and $\sigma'=s_{j-1}\sigma$.
Then
\begin{equation}
c_{j-1} c_{\sigma'} = c_{\sigma} + q c_{\sigma''},
\end{equation}
where $\sigma''$ is obtained from $\sigma'$ by exchanging $j-1$ with the
leftmost value smaller than $j-1$ to its right.

Moreover, the code of $\sigma''$ is obtained from the code of $\sigma'$ by
decrementing its $i$-th component, and $\sigma''$ still avoids the pattern 312.
\end{lemma}

\Proof
Let us first list all permutations $\tau$ with one inversion less than
$\sigma'$ such that $s_{j-1}\tau<\tau$, that is, such that $j-1$ is to the
right of $j$ in $\tau$.
Since $j-1$ is to the left of $j$ in $\sigma'$, we need to move one past the
other. Exchanging $j$ with an element to the left of $j-1$ in $\sigma'$ will
create inversions since all values to the left of $j-1$ in $\sigma'$ are
smaller that $j-1$ (proved above when justiyfing that $\sigma$ has a recoil in
$j-1$). So $j-1$ has to move to the right of $j$.

Now, we know that there are values smaller than $j-1$ to its right, we even
know that there are $a_i-1$ such values. Since they cannot be between $j-1$
and $j$, as  otherwise $\sigma$ would not avoid the 312 pattern, they have to be to
the right of $j$ in $\sigma'$. 
Since $\sigma'$ avoids the pattern $312$, all those values have to be
in decreasing order so that there is only one such value that can be exchanged
with $j-1$ and delete only one inversion: the first one.
So $\sigma''$ exists and is unique and is defined as in the statement.

The computation of the code of $\sigma''$ is straightforward and since we had
in $\sigma$ that $a_{i-1}\leq a_i-1$, we have here that $a''_i=a_i-2$ so that
$a''_i\geq a_{i-1}-1$, and $a''$ is again the code of a $312$-avoiding
permutation.
\qed

Now,
\begin{equation}
	\trg(c_{j-1}c_{\sigma'}) = \trg(c_{\sigma'}c_{j-1})=(q+1)\trg c_{\sigma'},
\end{equation}
whence the linear relation
\begin{equation}
(q+1)\trg c_{\sigma'} = \trg c_{\sigma} + q \trg c_{\sigma''}.
\end{equation}
Hence, all the local linear relations of the first type are induced by products
of the Kazhdan-Lusztig basis. The second type is dealt with similarly. \qed

\section{Annex: tables}

Here are tables of the $\tilde H_\mu$ with multi $t$ parameters up to $n=6$
with at least three rows and two columns.

\label{tables}
{\footnotesize

\begin{equation}
\widetilde{H}_{211}=
s_{4}+( q + t_{1} + t_{2} ) s_{31}+( q t_{1} + t_{2} ) s_{22}+( q t_{1} + q t_{2} + t_{1} t_{2} ) s_{211}+ q t_{1} t_{2}  s_{1111}
\end{equation}

\begin{equation}
\begin{split}
\widetilde{H}_{311}&=
s_{5}+( q^{2} + q + t_{1} + t_{2} ) s_{41}+( q^{2} t_{1} + q^{2} + q t_{1} + q t_{2} + t_{2} ) s_{32}\\
&+( q^{3} + q^{2} t_{1} + q^{2} t_{2} + q t_{1} + q t_{2} + t_{1} t_{2} ) s_{311}\\
&+( q^{3} t_{1} + q^{2} t_{1} + q^{2} t_{2} + q t_{1} t_{2} + q t_{2} ) s_{221}+( q^{3} t_{1} + q^{3} t_{2} + q^{2} t_{1} t_{2} + q t_{1} t_{2} ) s_{2111}\\
&+ q^{3} t_{1} t_{2}  s_{11111}
\end{split}
\end{equation}

\begin{equation}
\begin{split}
\widetilde{H}_{221}&=
s_{5}+( q t_{1} + q + t_{1} + t_{2} ) s_{41}+( q t_{1}^{2} + q^{2} + q t_{1} + t_{1} t_{2} + t_{2} ) s_{32}\\
&+( q^{2} t_{1} + q t_{1}^{2} + q t_{1} t_{2} + q t_{1} + q t_{2} + t_{1} t_{2} ) s_{311}\\
&+( q^{2} t_{1}^{2} + q^{2} t_{1} + q t_{1} t_{2} + t_{1}^{2} t_{2} + q t_{2} ) s_{221}+( q^{2} t_{1}^{2} + q^{2} t_{1} t_{2} + q t_{1}^{2} t_{2} + q t_{1} t_{2} ) s_{2111}\\
&+ q^{2} t_{1}^{2} t_{2}  s_{11111}
\end{split}
\end{equation}

\begin{equation}
\begin{split}
\widetilde{H}_{2111}&=
s_{5}+( q + t_{1} + t_{2} + t_{3} ) s_{41}+( q t_{1} + q t_{2} + t_{1} t_{3} + t_{2} + t_{3} ) s_{32}\\
&+( q t_{1} + q t_{2} + t_{1} t_{2} + q t_{3} + t_{1} t_{3} + t_{2} t_{3} ) s_{311}\\
&+( q t_{1} t_{2} + q t_{1} t_{3} + q t_{2} + t_{1} t_{3} + t_{2} t_{3} ) s_{221}+( q t_{1} t_{2} + q t_{1} t_{3} + q t_{2} t_{3} + t_{1} t_{2} t_{3} ) s_{2111}\\
&+ q t_{1} t_{2} t_{3}  s_{11111}
\end{split}
\end{equation}

\begin{equation}
\begin {split}
\widetilde{H}_{411} & =
s_{6}+( q^{3} + q^{2} + q + t_{1} + t_{2} ) s_{51}+( q^{4} + q^{3} t_{1} + q^{3} + q^{2} t_{1} + q^{2} t_{2} + q^{2} + q t_{1} + q t_{2} + t_{2} ) s_{42}\\
&+( q^{5} + q^{4} + q^{3} t_{1} + q^{3} t_{2} + q^{3} + q^{2} t_{1} + q^{2} t_{2} + q t_{1} + q t_{2} + t_{1} t_{2} ) s_{411}+( q^{4} t_{1} + q^{3} + q^{2} t_{1} + q^{2} t_{2} + q t_{2} ) s_{33}\\
&+( q^{5} t_{1} + q^{5} + 2 q^{4} t_{1} + q^{4} t_{2} + q^{4} + 2 q^{3} t_{1} + 2 q^{3} t_{2} + q^{2} t_{1} t_{2} + q^{2} t_{1} + 2 q^{2} t_{2} + q t_{1} t_{2} + q t_{2} ) s_{321}\\
&+( q^{6} + q^{5} t_{1} + q^{5} t_{2} + q^{4} t_{1} + q^{4} t_{2} + q^{3} t_{1} t_{2} + q^{3} t_{1} + q^{3} t_{2} + q^{2} t_{1} t_{2} + q t_{1} t_{2} ) s_{3111}\\
&+( q^{5} t_{1} + q^{4} t_{1} + q^{4} t_{2} + q^{3} t_{1} t_{2} + q^{2} t_{2} ) s_{222}+( q^{6} t_{1} + q^{5} t_{1} + q^{5} t_{2} + q^{4} t_{1} t_{2} + q^{4} t_{1} + q^{4} t_{2} + q^{3} t_{1} t_{2} + q^{3} t_{2} + q^{2} t_{1} t_{2} ) s_{2211}\\
&+( q^{6} t_{1} + q^{6} t_{2} + q^{5} t_{1} t_{2} + q^{4} t_{1} t_{2} + q^{3} t_{1} t_{2} ) s_{21111}+ q^{6} t_{1} t_{2}  s_{111111}
\end{split}
\end{equation}

\begin{equation}
\begin {split}
\widetilde{H}_{321} & =
s_{6}+( q^{2} + q t_{1} + q + t_{1} + t_{2} ) s_{51}+( q^{3} + 2 q^{2} t_{1} + q t_{1}^{2} + q^{2} + q t_{1} + q t_{2} + t_{1} t_{2} + t_{2} ) s_{42}\\
&+( q^{3} t_{1} + q^{3} + 2 q^{2} t_{1} + q t_{1}^{2} + q^{2} t_{2} + q t_{1} t_{2} + q t_{1} + q t_{2} + t_{1} t_{2} ) s_{411}+( q^{2} t_{1}^{2} + q^{3} + q^{2} t_{1} + q t_{2} + t_{1} t_{2} ) s_{33}\\
&+( q^{3} t_{1}^{2} + q^{4} + 3 q^{3} t_{1} + 2 q^{2} t_{1}^{2} + q^{2} t_{1} t_{2} + q^{2} t_{1} + 2 q^{2} t_{2} + 3 q t_{1} t_{2} + t_{1}^{2} t_{2} + q t_{2} ) s_{321}\\
&+( q^{4} t_{1} + q^{3} t_{1}^{2} + q^{3} t_{1} t_{2} + q^{3} t_{1} + q^{2} t_{1}^{2} + q^{3} t_{2} + 2 q^{2} t_{1} t_{2} + q t_{1}^{2} t_{2} + q t_{1} t_{2} ) s_{3111}\\
&+( q^{4} t_{1} + q^{3} t_{1}^{2} + q^{2} t_{1} t_{2} + q t_{1}^{2} t_{2} + q^{2} t_{2} ) s_{222}+( q^{4} t_{1}^{2} + q^{4} t_{1} + q^{3} t_{1}^{2} + q^{3} t_{1} t_{2} + q^{2} t_{1}^{2} t_{2} + q^{3} t_{2} + 2 q^{2} t_{1} t_{2} + q t_{1}^{2} t_{2} ) s_{2211}\\
&+( q^{4} t_{1}^{2} + q^{4} t_{1} t_{2} + q^{3} t_{1}^{2} t_{2} + q^{3} t_{1} t_{2} + q^{2} t_{1}^{2} t_{2} ) s_{21111}+ q^{4} t_{1}^{2} t_{2}  s_{111111}
\end{split}
\end{equation}

\begin{equation}
\begin {split}
\widetilde{H}_{3111} & =
s_{6}+( q^{2} + q + t_{1} + t_{2} + t_{3} ) s_{51}+( q^{2} t_{1} + q^{2} t_{2} + q^{2} + q t_{1} + q t_{2} + q t_{3} + t_{1} t_{3} + t_{2} + t_{3} ) s_{42}\\
&+( q^{3} + q^{2} t_{1} + q^{2} t_{2} + q^{2} t_{3} + q t_{1} + q t_{2} + t_{1} t_{2} + q t_{3} + t_{1} t_{3} + t_{2} t_{3} ) s_{411}+( q^{2} t_{1} + q^{2} t_{2} + q t_{1} t_{3} + q t_{2} + t_{3} ) s_{33}\\
&+( q^{3} t_{1} + q^{3} t_{2} + q^{2} t_{1} t_{2} + q^{2} t_{1} t_{3} + q^{2} t_{1} + 2 q^{2} t_{2} + q t_{1} t_{2} + q^{2} t_{3} + 2 q t_{1} t_{3} + q t_{2} t_{3} + q t_{2} + q t_{3} + t_{1} t_{3} + t_{2} t_{3} ) s_{321}\\
&+( q^{3} t_{1} + q^{3} t_{2} + q^{2} t_{1} t_{2} + q^{3} t_{3} + q^{2} t_{1} t_{3} + q^{2} t_{2} t_{3} + q t_{1} t_{2} + q t_{1} t_{3} + q t_{2} t_{3} + t_{1} t_{2} t_{3} ) s_{3111}\\
&+( q^{3} t_{1} t_{2} + q^{2} t_{1} t_{3} + q^{2} t_{2} + q t_{1} t_{3} + q t_{2} t_{3} ) s_{222}\\
&+( q^{3} t_{1} t_{2} + q^{3} t_{1} t_{3} + q^{3} t_{2} + q^{2} t_{1} t_{2} + q^{2} t_{1} t_{3} + q^{2} t_{2} t_{3} + q t_{1} t_{2} t_{3} + q t_{1} t_{3} + q t_{2} t_{3} ) s_{2211}\\
&+( q^{3} t_{1} t_{2} + q^{3} t_{1} t_{3} + q^{3} t_{2} t_{3} + q^{2} t_{1} t_{2} t_{3} + q t_{1} t_{2} t_{3} ) s_{21111}+ q^{3} t_{1} t_{2} t_{3}  s_{111111}
\end{split}
\end{equation}

\begin{equation}
\begin {split}
\widetilde{H}_{222} & =
s_{6}+( q t_{1} + q t_{2} + q + t_{1} + t_{2} ) s_{51}+( q^{2} t_{1} + q^{2} t_{2} + q t_{1} t_{2} + q^{2} + q t_{1} + t_{1}^{2} + q t_{2} + t_{1} t_{2} + t_{2}^{2} ) s_{42}\\
&+( q^{2} t_{1} t_{2} + q^{2} t_{1} + q t_{1}^{2} + q^{2} t_{2} + 2 q t_{1} t_{2} + q t_{2}^{2} + q t_{1} + q t_{2} + t_{1} t_{2} ) s_{411}+( q^{3} + q t_{1}^{2} + q t_{1} t_{2} + q t_{2}^{2} + t_{1} t_{2} ) s_{33}\\
&+( q^{3} t_{1} + q^{2} t_{1}^{2} + q^{3} t_{2} + 2 q^{2} t_{1} t_{2} + q t_{1}^{2} t_{2} + q^{2} t_{2}^{2} + q t_{1} t_{2}^{2} + q^{2} t_{1} + q t_{1}^{2} + q^{2} t_{2} + 2 q t_{1} t_{2} + t_{1}^{2} t_{2} + q t_{2}^{2} + t_{1} t_{2}^{2} ) s_{321}\\
&+( q^{3} t_{1} t_{2} + q^{2} t_{1}^{2} t_{2} + q^{2} t_{1} t_{2}^{2} + q^{2} t_{1}^{2} + 2 q^{2} t_{1} t_{2} + q t_{1}^{2} t_{2} + q^{2} t_{2}^{2} + q t_{1} t_{2}^{2} + q t_{1} t_{2} ) s_{3111}\\
&+( q^{3} t_{1} t_{2} + q^{2} t_{1}^{2} + q^{2} t_{1} t_{2} + q^{2} t_{2}^{2} + t_{1}^{2} t_{2}^{2} ) s_{222}\\
&+( q^{3} t_{1}^{2} + q^{3} t_{1} t_{2} + q^{2} t_{1}^{2} t_{2} + q^{3} t_{2}^{2} + q^{2} t_{1} t_{2}^{2} + q t_{1}^{2} t_{2}^{2} + q^{2} t_{1} t_{2} + q t_{1}^{2} t_{2} + q t_{1} t_{2}^{2} ) s_{2211}\\
&+( q^{3} t_{1}^{2} t_{2} + q^{3} t_{1} t_{2}^{2} + q^{2} t_{1}^{2} t_{2}^{2} + q^{2} t_{1}^{2} t_{2} + q^{2} t_{1} t_{2}^{2} ) s_{21111}+ q^{3} t_{1}^{2} t_{2}^{2}  s_{111111}
\end{split}
\end{equation}

\begin{equation}
\begin {split}
\widetilde{H}_{2211} & =
s_{6}+( q t_{1} + q + t_{1} + t_{2} + t_{3} ) s_{51}+( q t_{1}^{2} + q t_{1} t_{2} + q^{2} + q t_{1} + q t_{2} + 2 t_{1} t_{3} + t_{2} + t_{3} ) s_{42}\\
&+( q^{2} t_{1} + q t_{1}^{2} + q t_{1} t_{2} + q t_{1} t_{3} + q t_{1} + q t_{2} + t_{1} t_{2} + q t_{3} + t_{1} t_{3} + t_{2} t_{3} ) s_{411}+( q^{2} t_{1} + q t_{1} t_{2} + t_{1}^{2} t_{3} + q t_{2} + t_{3} ) s_{33}\\
&+( q^{2} t_{1}^{2} + q^{2} t_{1} t_{2} + q t_{1}^{2} t_{2} + q t_{1}^{2} t_{3} + q^{2} t_{1} + q^{2} t_{2} + 2 q t_{1} t_{2} + 2 q t_{1} t_{3} + t_{1}^{2} t_{3} + t_{1} t_{2} t_{3} + q t_{2} + q t_{3} + t_{1} t_{3} + t_{2} t_{3} ) s_{321}\\
&+( q^{2} t_{1}^{2} + q^{2} t_{1} t_{2} + q t_{1}^{2} t_{2} + q^{2} t_{1} t_{3} + q t_{1}^{2} t_{3} + q t_{1} t_{2} t_{3} + q t_{1} t_{2} + q t_{1} t_{3} + q t_{2} t_{3} + t_{1} t_{2} t_{3} ) s_{3111}\\
&+( q^{2} t_{1}^{2} t_{2} + q t_{1}^{2} t_{3} + q^{2} t_{2} + q t_{1} t_{3} + t_{1} t_{2} t_{3} ) s_{222}\\
&+( q^{2} t_{1}^{2} t_{2} + q^{2} t_{1}^{2} t_{3} + 2 q^{2} t_{1} t_{2} + q t_{1}^{2} t_{3} + q t_{1} t_{2} t_{3} + t_{1}^{2} t_{2} t_{3} + q t_{1} t_{3} + q t_{2} t_{3} ) s_{2211}\\
&+( q^{2} t_{1}^{2} t_{2} + q^{2} t_{1}^{2} t_{3} + q^{2} t_{1} t_{2} t_{3} + q t_{1}^{2} t_{2} t_{3} + q t_{1} t_{2} t_{3} ) s_{21111}+ q^{2} t_{1}^{2} t_{2} t_{3}  s_{111111}
\end{split}
\end{equation}

\begin{equation}
\begin {split}
\widetilde{H}_{21111} & =
s_{6}+( q + t_{1} + t_{2} + t_{3} + t_{4} ) s_{51}+( q t_{1} + q t_{2} + q t_{3} + t_{1} t_{3} + t_{1} t_{4} + t_{2} t_{4} + t_{2} + t_{3} + t_{4} ) s_{42}\\
&+( q t_{1} + q t_{2} + t_{1} t_{2} + q t_{3} + t_{1} t_{3} + t_{2} t_{3} + q t_{4} + t_{1} t_{4} + t_{2} t_{4} + t_{3} t_{4} ) s_{411}+( q t_{1} t_{3} + q t_{2} + t_{1} t_{4} + t_{2} t_{4} + t_{3} ) s_{33}\\
&+( q t_{1} t_{2} + 2 q t_{1} t_{3} + q t_{2} t_{3} + q t_{1} t_{4} + q t_{2} t_{4} + t_{1} t_{2} t_{4} + t_{1} t_{3} t_{4} + q t_{2} + q t_{3} + t_{1} t_{3} + t_{2} t_{3} + t_{1} t_{4} + 2 t_{2} t_{4} + t_{3} t_{4} ) s_{321}\\
&+( q t_{1} t_{2} + q t_{1} t_{3} + q t_{2} t_{3} + t_{1} t_{2} t_{3} + q t_{1} t_{4} + q t_{2} t_{4} + t_{1} t_{2} t_{4} + q t_{3} t_{4} + t_{1} t_{3} t_{4} + t_{2} t_{3} t_{4} ) s_{3111}\\
&+( q t_{1} t_{2} t_{4} + q t_{1} t_{3} + q t_{2} t_{3} + t_{1} t_{3} t_{4} + t_{2} t_{4} ) s_{222}\\
&+( q t_{1} t_{2} t_{3} + q t_{1} t_{2} t_{4} + q t_{1} t_{3} t_{4} + q t_{1} t_{3} + q t_{2} t_{3} + q t_{2} t_{4} + t_{1} t_{2} t_{4} + t_{1} t_{3} t_{4} + t_{2} t_{3} t_{4} ) s_{2211}\\
&+( q t_{1} t_{2} t_{3} + q t_{1} t_{2} t_{4} + q t_{1} t_{3} t_{4} + q t_{2} t_{3} t_{4} + t_{1} t_{2} t_{3} t_{4} ) s_{21111}+ q t_{1} t_{2} t_{3} t_{4}  s_{111111}
\end{split}
\end{equation}

}


\end{document}